\documentclass[11pt]{article}
\usepackage[margin=1.5in]{geometry}
\usepackage{amssymb}
\usepackage{enumitem}
\usepackage{amsmath}
\usepackage{amsthm}
\usepackage{graphicx}
\usepackage{caption}
\usepackage{siunitx}
\usepackage{mathrsfs}
\usepackage{quiver}
\usepackage{tikz}
\usepackage{tikz-cd}
\usepackage{adjustbox}
\usepackage{hyperref}
\usepackage{bbold}
\usepackage[affil-it]{authblk}
\usepackage{qworld}
\usepackage{qworld-preview}

\DeclareMathOperator{\Hom}{Hom}
\DeclareMathOperator{\Cone}{Cone}

\DeclareMathOperator{\End}{End}
\DeclareMathOperator{\id}{id}

\DeclareMathOperator{\Ch}{Ch}
\DeclareMathOperator{\SBim}{SBim}
\DeclareMathOperator{\JM}{JM}
\DeclareMathOperator{\FT}{FT}
\DeclareMathOperator{\tw}{tw}

\DeclareMathOperator{\pTr}{pTr}

\newtheorem{thm}{Theorem}[section]
\newtheorem{proposition}[thm]{Proposition}

\newtheorem{lemma}[thm]{Lemma}
\newtheorem{corollary}[thm]{Corollary}

\theoremstyle{definition}
\newtheorem{definition}{Definition}[section]

\theoremstyle{remark}
\newtheorem*{remark}{Remark}

\numberwithin{equation}{section}

\usetikzlibrary{calc, shapes, backgrounds, positioning}
\title{Recursive relations for categorified Young symmetrizers
in type $\mathbb{B}$}
\author{Nur \c{C}ataltepe\thanks{\texttt{cataltepe.a@northeastern.edu}}}
\affil{Department of Mathematics, Northeastern University}
\date{November 20, 2026}
\begin{document}
\maketitle
\begin{abstract}
\noindent The category of bounded-above chain complexes
of Soergel bimodules of any finite type $\mathbb{X}_n$
contains a symmetrizer complex $\mathbf{P}_{\mathbb{X}_n}$ which
kills non-identity Soergel bimodules up to homotopy. It has been shown in
previous work that, in type $\mathbb{A}$,
$\mathbf{P}_{\mathbb{A}_{n+1}}$ can be obtained as a
semi-infinite
extension between copies of
$\mathbf{P}_{\mathbb{A}_n}$ and its product with a
Jucys-Murphy complex. We prove
an analogous recursion holds between the symmetrizers for
type $\mathbb{B}_n$ and $\mathbb{B}_{n+1}$. Following the method of the result in type $\mathbb{A}$,
we then use this to prove that the cohomology ring of the
endomorphism algebra of
$\mathbf{P}_{\mathbb{B}_n}$ is a polynomial ring.
\end{abstract}
\section{Introduction}
\noindent Categorified Young symmetrizers
in the category of chain complexes of
Soergel bimodules in type $\mathbb{A}$
have proven useful to calculating the
triply-graded homology of
links~\cite{hogancamp_categorified_2018, abel_categorified_2016,
elias_computation_2019} and
play an important role in showing the Rouquier complexes corresponding to
full-twist braids are categorically diagonalizable~\cite{elias_categorical_2025}.
Central to this are the recursive relations 
of the following form between finite projectors (which kill
non-identity Soergel bimodules within their corresponding parabolic subgroup)
and Jucys-Murphy Rouquier complexes (omitting some grading shifts which are all
parity in cohomology):
\[K_{\mathbb{A}_n}\JM_{n+1} \simeq K_{\mathbb{A}_{n+1}} \to K_{\mathbb{A}_{n}}\]
This recursion follows from a recursion between the corresponding Young symmetrizers
$\mathbf{P}_{\mathbb{A}_n}$
and $\mathbf{P}_{\mathbb{A}_{n+1}}$ (defined in Theorem~\ref{thm:sbimidempotent})
where, once again
omitting some even-cohomological grading shifts, we can obtain
$\mathbf{P}_{\mathbb{A}_{n+1}}$ as a semi-infinite extension
\[\begin{tikzcd}[ampersand replacement=\&,cramped]
	{\mathbf{P}_{\mathbb{A}_{n}}} \& {\mathbf{P}_{\mathbb{A}_{n}}\JM_{n+1}} \\
	{\mathbf{P}_{\mathbb{A}_n}} \& {\mathbf{P}_{\mathbb{A}_n}\JM_{n+1}} \\
	{\mathbf{P}_{\mathbb{A}_n}} \& {\mathbf{P}_{\mathbb{A}_n}\JM_{n+1}} \\
	\dots \& \dots
	\arrow[from=1-1, to=1-2]
	\arrow[from=1-1, to=2-2]
	\arrow[from=2-1, to=2-2]
	\arrow[from=2-1, to=3-2]
	\arrow[from=3-1, to=3-2]
\end{tikzcd}\]
This construction in turn
is a result of $\mathbf{P}_{\mathbb{A}_n}\JM_{n+1}$ being categorically
diagonalizable with two eigenmaps, which is motivated by the fact that
$\Ch^-(\SBim_{\mathbb{A}_n})\mathbf{P}_{\mathbb{A}_n}$ treated as a categorical left module of
$\SBim_{\mathbb{A}_n}$, has only two left cells.
\\ \\
Now, given the Young symmetrizer
$\mathbf{P}_{\mathbb{B}_n}$ in type $\mathbb{B}_n$, one can observe that
$\Ch^-(\SBim_{\mathbb{B}_n})\mathbf{P}_{\mathbb{B}_n}$ viewed as a categorical left module of
$\SBim_{\mathbb{B}_n}$ also has two left cells. This motivates us to
prove a version of the above recursion for type $\mathbb{B}$ of
the exact same form (with different grading shifts and an appropriate
choice of generalized Jucys-Murphy braid), which we inductively do
in the statements in
Section~\ref{sec:recursiverelations}, following methods in~\cite{elias_categorical_2017, hogancamp_categorified_2018, abel_categorified_2016,elias_computation_2019}.
We can then use our recursive relations to prove Theorem~\ref{thm:ringstructure},
which shows that $H^\bullet\End(\mathbf{P}_{\mathbb{B}_n})$ is isomorphic to
a bigraded polynomial ring in $n$ variables.
\\ \\
Our results represent a first step towards diagonalizing full-twist complexes
in type $\mathbb{B}$.
\subsection*{Acknowledgements}
This problem was suggested by the author's Ph.D. advisor
Matthew Hogancamp, and the progress
towards these results
was greatly aided by conversations with
him.
\\ \\
\textbf{AI disclosure:}
Generative AI was not used in the process of deriving
these results or writing this paper.
\section{Preliminaries}
We are working in the $dg$-category of bounded-above cochain complexes
over the category
$\SBim_{\mathbb{B}_n}$ of Soergel bimodules of type
$\mathbb{B}_n$, denoted
$\Ch^-(\SBim_{\mathbb{B}_n})$, and will occasionally refer to its homotopy category
$\mathcal{K}^-(\SBim_{\mathbb{B}_n})$. We refer the
reader to the textbook~\cite{elias_introduction_2020} for an introduction
to Soergel bimodules and the diagrammatics of morphisms between them.
Throughout this paper, we will assume our
bimodules to be defined using the signed permutation realization of
the Coxeter system of type $\mathbb{B}_n$, though our results hold for other
realizations with minimal modifications\footnote{In fact, the only point where
we rely on a specific fact about a realization is
Lemma~\ref{lem:subcomplexeliminations}, and this fact continues to hold in
the geometric realization. It should, however, be noted that the realization
we use in this paper allows us to work over
$\mathbb{Q}$.}.
Objects in this category are bigraded, with the cohomological shift functor
(and grading)
denoted by $T$ and the bimodule grading shift functor (and grading) denoted by $Q$.
When we are making statements about general categories of Soergel bimodules, we
may replace $\mathbb{B}_n$ with $\mathbb{X}_n$ for some Dynkin diagram of type $\mathbb{X}_n$ or
$(W,S)$ for a Coxeter system consisting of a group $W$ and set of simple reflections
$S$.
\subsection{How we write extensions of complexes}
Throughout this paper, we will be
frequently dealing with complexes realized as twists
of direct sums of other complexes, and will be adopting certain conventions to
make this easier. We will write a cochain complex that is a twist of a direct
sum of two complexes as 
\[X\xrightarrow{\alpha} Y\ \mathrm{or}\ \tw_{\alpha}(X\oplus Y)\]
where $X$ and $Y$ are cochain complexes with respective differentials
$d_X$ and $d_Y$.
Unless otherwise noted, $\alpha$ is a bimodule map of $Q$-degree $0$ and
$T$-degree $1$, such that
$\left(\begin{smallmatrix} d_X & 0 \\ \alpha & d_Y \end{smallmatrix}\right)^2=0$
(using matrix notation to write the differentials
and $\alpha$ as endomorphisms of $X\oplus Y$). This is equivalent to
\[\Cone(TX\xrightarrow{\alpha} Y)\]
where we have abused notation (as we similarly will later) to equate the chain map
$\alpha:TX\to Y$ with the degree-1 map $\alpha:X\to Y$.
\\ \\
When writing a twist $Z$ of a direct sum of a (possibly infinite) family of complexes
$\bigoplus_{i=1}^\infty X_i$ such that $d_Z$ only includes a degree-$1$ morphism
$X_i\to X_j$ if $i \leq j$, with the component being the original differential
$d_{X_i}$ itself if $i=j$, we will write the resultant complex as follows
\[\begin{tikzcd}[ampersand replacement=\&,cramped]
	{X_1} \& {X_2} \& \\
	\& {X_3} \& {X_4} \\
	\& \dots \& \dots
	\arrow[from=1-1, to=1-2]
	\arrow[from=1-1, to=2-2]
	\arrow[from=1-2, to=2-2]
	\arrow[from=2-2, to=2-3]
	\arrow[from=2-3, to=3-3]
\end{tikzcd}\]
where the arrows drawn may be a strict subset of the actual components of $d_Z$.
Arrows will be labeled only if they are components of $d_Z$ with
known values. We will often
write complexes in this way without comment, where the numbering scheme of the
direct sum can be inferred from the directedness of the graph created by the arrows.
\\ \\
In some cases where some large contraction of intermediate terms in such a complex
has occurred, we may write complexes of the form
\begin{align*}
(X &\to Y )\\
	&{\big\downarrow}\\
(Z &\to V \to W\to \dots)
\end{align*}
where the vertical arrow between the parenthesized regions is not labeled and
merely indicates the probable presence of an unknown, degree-$1$ arrow in the
differential of the total complex that is between the complexes represented by the
two parenthesized regions.
\subsection{New grading conventions}
For
notational convenience, we now introduce new grading shifts
$t:=T^2Q^{-2}$ and $q:=Q^{2}$, and will be using these throughout the rest of
this paper. Note that this means $T=(tq)^{1/2}$ and $Q=q^{1/2}$.
\subsection{Hecke category shorthands}
We will be using the graphical calculus for maps between Soergel bimodules developed
by Elias and Williamson~\cite{elias_introduction_2020}, and will in particular be
extensively using the following shorthand to indicate the idempotent in the
endomorphism ring of the Bott-Samelson bimodule
$\mathrm{BS}(s_{i+1}s_is_{i+1})$ corresponding to the
indecomposable Soergel bimodule $B_{s_{i+1}s_{i}s_{i+1}}$ in
type $\mathbb{B}_n$, where $1 \leq i \leq n-1$:
\[\q[scale=0.5]{
		\qwire[color=pink]\qwire[color=black]\qwire[color=pink]\n
		\qbox[n=3,s=3,N={},S={},hlen=3,id=idempbox]\n
		\qwire[color=pink]\qwire[color=black]\qwire[color=pink]
		\qwire[cod=s-1-idempbox,color=pink]
		\qwire[cod=s-2-idempbox,color=black]
		\qwire[cod=s-3-idempbox,color=pink]
		\qwire[dom=n-1-idempbox,color=pink]
		\qwire[dom=n-2-idempbox,color=black]
		\qwire[dom=n-3-idempbox,color=pink]
}:=\q[scale=0.5]{
	\qwire[color=pink]\qwire[color=black]\qwire[color=pink]\n
	\qwire[color=pink]\qwire[color=black]\qwire[color=pink]\n
	\qwire[color=pink]\qwire[color=black]\qwire[color=pink]
}-\frac{1}{\partial_{i+1}(\alpha_i)}
\q[scale=0.5]{
    \qwire[color=pink, id=upwire1] \qunit[color=black, frame color=black]
    \qwire[color=pink, id=upwire2]\n
    \qbox[n=3,s=3,N={},S={},hlen=3,id=idempbox,frame color=none]\n
    \qwire[color=pink, id=downwire1] \qcounit[color=black, frame color=black]
    \qwire[color=pink, id=downwire2]
    \qwire[color=pink, dom=N-downwire1, cod=s-2-idempbox]
    \qwire[color=pink, dom=N-downwire2, cod=s-2-idempbox]
    \qwire[color=pink, dom=s-2-idempbox, cod=n-2-idempbox]
    \qwire[color=pink, cod=S-upwire1, dom=n-2-idempbox]
    \qwire[color=pink, cod=S-upwire2, dom=n-2-idempbox]
}\]
In this figure, the black strands correspond to $s_i$ and the pink strands
to $s_{i+1}$. We will be representing morphisms between Soergel bimodules
using diagrams representing morphisms between the corresponding Bott-Samelsons,
composed with the relevant projectors.
\subsection{Induction and partial trace functors in type $\mathbb{B}$}
There is a natural induction/``inclusion'' functor $\SBim_{\mathbb{B}_n}\to \SBim_{\mathbb{B}_{n+1}}$.
\begin{definition}
	\label{def:inductionfunctor}
	We denote by $I_n$ the additive \textbf{induction} functor
	$\SBim_{\mathbb{B}_n}\to \SBim_{\mathbb{B}_{n+1}}$ such that
	\begin{itemize}
		\item $I_n(B_x)$ for an
			indecomposable Soergel bimodule corresponding
			to an element $x$ of the Coxeter group is
			$B_x\in\SBim_{\mathbb{B}_{n+1}}$.
		\item $I_n(f)$ for a morphism $f:B_x\to B_y$ is
			that same $f$ but now living in a bimodule over
			the symmetric algebra of one rank higher.
	\end{itemize}
	Thus, $I_n$ is faithful on hom spaces.
As an additive functor, $I_n$ extends to a $dg$-functor
$\Ch^-(\SBim_{\mathbb{B}_{n}})\to \Ch^-(\SBim_{\mathbb{B}_{n+1}})$.
\end{definition}
\noindent We will also be making use of the right adjoint of $I_n$, which can most
naturally be defined over the category of cochain complexes.
\begin{definition}
	\label{def:partialtrace}
	We denote by $\pTr_n$ the \textbf{partial trace}
	$\Ch^-(\SBim_{\mathbb{B}_{n+1}})\to \Ch^-(\SBim_{\mathbb{B}_n})$, which is the
	right adjoint of $I_n$.
\end{definition}
\noindent In~\cite{cataltepe_partial_2026} we observe that
$\pTr_n$
is well-behaved on Soergel bimodules in type $\mathbb{B}$.
\begin{lemma}[\cite{cataltepe_partial_2026}]
	\label{lem:partialtrace}
	The object $\pTr_{n}(X)$ for a Soergel
	bimodule $X$ is naturally homotopy
	equivalent to $\bigoplus_{i=0}^\infty q^i Y$ for some fixed
	Soergel bimodule $Y\in \SBim_{\mathbb{B}_n}$, with the infinite $q$-graded
	sum representing the extra variable added to $\mathbb{1}$ in type
	$\mathbb{B}_{n+1}$.
\end{lemma}
\subsection{Unital idempotents in $\Ch^-(\SBim_{\mathbb{X}_n})$}
We are concerned with a particular class of cochain complexes of
Soergel bimodules.
\begin{thm}[\cite{hogancamp_constructing_2020}]
	\label{thm:sbimidempotent}
	For any finite type Coxeter system $\mathbb{X}_n$, there exists a unique
	(up to homotopy equivalence)
	pair $(\mathbf{P}_{\mathbb{X}_n},\eta_{\mathbb{X}_n})$ for
	$\mathbf{P}_{\mathbb{X}_n} \in \Ch^-(\SBim_{\mathbb{X}_n})$ and
	a cochain map
	$\eta_{\mathbb{X}_n}:\mathbb{1}\to \mathbf{P}_{\mathbb{X}_n}$ such that
	\begin{enumerate}[label=(\arabic*)]
		\item $\mathbf{P}_{\mathbb{X}_n}$ is in the Drinfeld center of
			$\mathcal{K}^-(\SBim_{\mathbb{X}_n})$
		\item $\mathbf{P}_{\mathbb{X}_n}B_{s}\simeq B_{s}
			\mathbf{P}_{\mathbb{X}_n}\simeq 0$ for any non-identity
			simple reflection $s$.
		\item $\Cone(\eta_{\mathbb{X}_n})$ is homotopy equivalent to
			a complex that does not contain any shifts of
			$\mathbb{1}$ as summands.
	\end{enumerate}	
\end{thm}
\noindent The conditions (2) and (3) ensure that $\mathbf{P}_{\mathbb{X}_n}
\simeq \mathbf{P}_{\mathbb{X}_n}\otimes \mathbf{P}_{\mathbb{X}_n}$, since
\begin{align*}\Cone(\mathbf{P}_{\mathbb{X}_n}
	\xrightarrow{\id_{\mathbf{P}_{\mathbb{X}_n}}
	\otimes \eta_{\mathbb{X}_n}}
	\mathbf{P}_{\mathbb{X}_n}\otimes \mathbf{P}_{\mathbb{X}_n})&\\=
	\mathbf{P}_{\mathbb{X}_n}\otimes \Cone(\eta_{\mathbb{X}_n}) &\simeq 0 \simeq
	\Cone(\eta_{\mathbb{X}_n})\otimes \mathbf{P}_{\mathbb{X}_n}\\ &=
	\Cone(\mathbf{P}_{\mathbb{X}_n}\xrightarrow{\eta_{\mathbb{X}_n}\otimes
\id_{\mathbf{P}_{\mathbb{X}_n}}}\mathbf{P}_{\mathbb{X}_n}
\otimes \mathbf{P}_{\mathbb{X}_n})\end{align*}
meaning the maps $\id_{\mathbf{P}_{\mathbb{X}_n}} \otimes
\eta_{\mathbb{X}_n}$ and $\eta_{\mathbb{X}_n}\otimes \id_{\mathbf{P}_{\mathbb{X}_n}}$ are homotopy equivalences.
This makes $(\mathbf{P}_{\mathbb{X}_n}, \eta_{\mathbb{X}_n})$ a unital idempotent in
$\Ch^-(\SBim_{\mathbb{X}_n})$, the definition of which we recall now.
\begin{definition}
	A \textbf{unital idempotent} $(P,\eta)$ in a $dg$-category is an object $P$
and a closed degree-zero morphism $\eta:\mathbb{1}\to P$ such that the maps
\[P\xrightarrow{\id_P\otimes \eta} P\otimes P\ \mathrm{and}\ 
P\xrightarrow{\eta\otimes \id_P} P\otimes P\] are
homotopy equivalences.
Equivalently, $P\otimes \Cone(\eta)\simeq \Cone(\eta)\otimes P\simeq  0$.
\end{definition}
\begin{remark}
	If the Coxeter system of type $\mathbb{X}_n$ is a parabolic subgroup of
	some other one of type $\mathbb{Y}_m$, versions of the properties of
	$\mathbf{P}_{\mathbb{X}_n}$ continue to apply to its image under the natural
	induction functor between the categories of cochain complexes
	of Soergel bimodules. That is, the image of $\mathbf{P}_{\mathbb{X}_n}$
	eliminates all elementary Bott-Samelsons belonging to the parabolic
	subgroup of type $\mathbb{X}_n$, and the cone of the image of
	$\eta_{\mathbb{X}_n}$ consists only of non-identity Soergel bimodules in the
	parabolic subgroup.
\end{remark}
\noindent A useful fact about categorical idempotents is the following.
\begin{thm}[\cite{hogancamp_constructing_2020, hogancamp_idempotents_2017}]
\label{thm:idempotentmorphismzero}
Let $(P,\eta)$ be a unital idempotent in a $dg$-category
$\mathcal{C}$, and let $M,N\in\mathcal{C}$ such that
\begin{itemize}
	\item $PM\simeq 0$ and 
	\item $PN\simeq N$.
\end{itemize}
\noindent Then, $\Hom_{\mathcal{C}}(M,N)\simeq 0$.
\end{thm}
\noindent This fact 
	gives a homotopy
	equivalence $\End(\mathbf{P}_{\mathbb{X}_n})\simeq \Hom(\mathbb{1},
	\mathbf{P}_{\mathbb{X}_n})$ given by $f\mapsto f\circ \eta_{\mathbb{X}_n}$, 
	since
	\[\Cone\left(\End(\mathbf{P}_{\mathbb{X}_n})
			\xrightarrow{\eta_{\mathbb{X}_n}^\star}
		\Hom(\mathbb{1},\mathbf{P}_{\mathbb{X}_n})\right)=
	\Hom(\Cone(\eta_{\mathbb{X}_n}), \mathbf{P}_{\mathbb{X}_n}),\]
	which is therefore contractible since $\mathbf{P}_{\mathbb{X}_n}
	\otimes \Cone(\eta_{\mathbb{X}_n})
\simeq 0$.
\subsection{Generalized Jucys-Murphy complexes}
We begin by recalling the definition of a Rouquier complex.
\begin{definition}
Given a Coxeter system $(W,S)$ and a simple reflection $s\in S$, the
\textbf{elementary Rouquier complex} $F(s)\in
\Ch^-(\SBim_{(W,S)})$ is given by
\[B_s \xrightarrow{\q[scale=0.4]{\qcounit[color=black]}} t^{1/2}\mathbb{1}\]
where $B_{s}$ is the elementary Bott-Samelson bimodule corresponding to
$s$. We denote
\[F(s)^{-1}:=t^{-1/2}\mathbb{1} 
\xrightarrow{\q[scale=0.4]{\qunit[color=black]}} B_{s}.\]
Given a braid word $\sigma=\sigma_{1}^{e_1}\dots \sigma_k^{e_k}$ in the
braid group corresponding to $(W,S)$, where $\sigma_i$ is the
positive braid lift of some $s_i\in S$ and $e_i \in \{1,-1\}$, we write
the \textbf{Rouquier complex} $F(\sigma)$ to be the tensor product
\[F(s_1)^{e_1}\dots F(s_k)^{e_k}.\]
We denote the inverse of this Rouquier complex by
\[F(\sigma)^{-1}:=F(s_k)^{-e_k}\dots F(s_1)^{-e_1}=F(\sigma^{-1}).\]
Rouquier complexes satisfy the braid relations in $W$ up to homotopy equivalence,
and $F(\sigma)F(\sigma^{-1})\simeq F(\sigma^{-1})F(\sigma)\simeq \mathbb{1}$.
\end{definition}
\noindent If $w$ is an element of $(W,S)$, we will denote the Rouquier
complex corresponding to the positive braid lift of any of its reduced expressions
(which is well-defined up to homotopy equivalence) by $F(w)$ for short.
\begin{definition}
Given a finite-type Coxeter system $W$ of type
$\mathbb{X}_n$ with longest element $w_0$, the
\textbf{full twist} $\FT(\mathbb{X}_n)$ or $\FT(W)$ is the Rouquier complex
$F(w_0)^2$, which is well-defined up to homotopy equivalence.
\end{definition}
\begin{definition}
\noindent Given a parabolic subgroup $W^I$ of a finite-type
Coxeter system $W$, the
\textbf{generalized Jucys-Murphy complex from $W^I$ to $W$} denoted
$\JM(W^I,W)$ (or $\JM(\mathbb{X}_n,\mathbb{Y}_m)$ if the types of $W^I$ and $W$ more
clearly indicate the nature of the complex) is the Rouquier complex
\[\JM(W^I,W):=
\FT(W)\FT(W^I)^{-1}\simeq \FT(W^I)^{-1}\FT(W)\] (as $\FT(W)$ is central).
\end{definition}
\begin{remark}
The full twist of $W$ is central in the braid group of $W$, and
the full twist of
$W^I$ is central in the braid subgroup given by the braid group of
$W^I$, so the complexes $\FT(W)$ and $\FT(W^I)$, and therefore
$\JM(W^I,W)$ all commute up to homotopy equivalence with any cochain
complex consisting of Soergel bimodules in $W^I$.
\end{remark}
\noindent The generalized Jucys-Murphy complex from type
$\mathbb{B}_n$ to $\mathbb{B}_{n+1}$ is given
by 
\begin{equation}\label{eqn:typebjm}\JM(\mathbb{B}_n, \mathbb{B}_{n+1})=
	F(s_{n+1}s_{n}\dots s_{2}s_1s_2
\dots s_{n}s_{n+1})^2.\end{equation}
\subsubsection{Eigenmaps and recursive relations in type $\mathbb{A}$}
The result of this paper adapts results in~\cite{elias_computation_2019} in type
$\mathbb{A}$
giving recursive relations between the categorified Young symmetrizers in
$\mathbb{A}_n$ and $\mathbb{A}_{n+1}$.
These symmetrizers are crucial to proving results about
categorical diagonalization in type $\mathbb{A}$, as well as calculating triply-graded
link homology for links. We hope the existence of similar recursive
relations in type $\mathbb{B}$ will allow us to adapt results about categorical
diagonalization to full twist complexes in type $\mathbb{B}$, where no connection to
knot theory exists.
\section{Recursive relations for categorified Young symmetrizers in type $\mathbb{B}$}
\label{sec:recursiverelations}
We now establish recursive relations between the symmetrizer complexes
in type $\mathbb{B}$. We will use $P_{\mathbb{B}_n}\in
\Ch^-(\SBim_{\mathbb{B}_{n+1}})$ to denote the
image of $\mathbf{P}_{\mathbb{B}_n}\in \Ch^-(\SBim_{\mathbb{B}_n})$
(as defined in Theorem~\ref{thm:sbimidempotent}) under the induction functor $I_n$.
\begin{lemma}
	\label{lem:eigenmaps}
	For each $n \geq 1$, there exist closed maps
	\[\alpha_n,\beta_n: P_{\mathbb{B}_n}\to P_{\mathbb{B}_n}
	\JM(\mathbb{B}_n,\mathbb{B}_{n+1})\]
	where $\deg(\beta_n)=q$ and $\deg(\alpha_n)=t^{2n+1}$, such that:
	\begin{itemize}
		\item $\Cone(\alpha_n)$ is homotopy equivalent to
			a complex whose terms are all shifts of
			non-identity Soergel bimodules.
		\item $\Cone(\beta_n)B_s\simeq B_s\Cone(\beta_n)\simeq 0$
			for all elementary Bott-Samelson bimodules
			$B_s\in \SBim_{\mathbb{B}_{n+1}}$.
	\end{itemize}
\end{lemma}
\noindent To prove this lemma,
 we compute a simplified form for the complex
 $P_{\mathbb{B}_n}\JM(\mathbb{B}_n,\mathbb{B}_{n+1})$ and then show that:
\begin{itemize}
	\item There exists a map $\alpha_n$ which is induced by the inclusion of
		the last term $t^{2n+1}P_{\mathbb{B}_n}$ into
		$P_{\mathbb{B}_n}\JM(\mathbb{B}_n,\mathbb{B}_{n+1})$.
	\item There exists a map $\beta_n:qP_{\mathbb{B}_n}
		\to P_{\mathbb{B}_n}\JM(\mathbb{B}_{n},\mathbb{B}_{n+1})$.
		with a component going to the first term in the latter
		complex of the form
		\[qP_{\mathbb{B}_n} \xrightarrow{\id_{P_{\mathbb{B}_n}} \otimes
		\q[scale=0.5]{\qunit[color=black]}} q^{1/2}P_{\mathbb{B}_n}
	B_{s_{n+1}}.\]
\end{itemize}
The desired properties of the mapping cones then follow.
The existence of these maps immediately affords the following property of
$P_{\mathbb{B}_n}\JM(\mathbb{B}_n,\mathbb{B}_{n+1})$:
\begin{corollary}
	$P_{\mathbb{B}_n}\JM(\mathbb{B}_n,\mathbb{B}_{n+1})$
	is categorically diagonalizable within
	the category $P_{\mathbb{B}_n}\Ch^-(\SBim_{\mathbb{B}_{n+1}})
	P_{\mathbb{B}_n}$ of complexes
	expressible as tensor products of complexes of the form
	$P_{\mathbb{B}_n}XP_{\mathbb{B}_n}$, where
	$X\in\Ch^-(\SBim_{\mathbb{B}_{n+1}})$. Its eigenmaps
	are given by the maps $\alpha_n,\beta_n$ in
	Lemma~\ref{lem:eigenmaps}.
\end{corollary}
\noindent We may then apply methods from categorical diagonalization to
give an expression of $\mathbf{P}_{\mathbb{B}_{n+1}}$ in terms of a
semi-infinite exitension between the cones of
these eigenmaps.
\begin{proposition}
	\label{prop:symexistence}
	Given $n\geq 1$, let $\alpha_n,\beta_n$ be as in
	Lemma~\ref{lem:eigenmaps}.
	Let $P_{\mathbb{B}_{n+1}}'$
	be the complex constructed as the following semi-infinite
	extension of complexes:
\begin{equation}\label{eq:nextstepconstruction}
	P_{\mathbb{B}_{n+1}}':=	t^{-(2n+1)}\left(
\begin{tikzcd}[ampersand replacement=\&,cramped]
	{(tq)^{-1/2}qP_{\mathbb{B}_n}} \& {P_{\mathbb{B}_n}\JM(\mathbb{B}_n,\mathbb{B}_{n+1})} \\
	{(tq)^{-1/2}q^2t^{-(2n+1)}P_{\mathbb{B}_n}} \& {qt^{-(2n+1)}P_{\mathbb{B}_n}\JM(\mathbb{B}_n,\mathbb{B}_{n+1})} \\
	{(tq)^{-1/2}
	q^3t^{-2(2n+1)}P_{\mathbb{B}_n}} \& {q^2t^{-2(2n+1)}P_{\mathbb{B}_n}\JM(\mathbb{B}_n,\mathbb{B}_{n+1})} \\
	\dots \& \dots
	\arrow["{\beta_n}", from=1-1, to=1-2]
	\arrow["{\alpha_n}", from=1-1, to=2-2]
	\arrow["{\beta_n}", from=2-1, to=2-2]
	\arrow["{\alpha_n}", from=2-1, to=3-2]
	\arrow["{\beta_n}", from=3-1, to=3-2]
	\arrow["{\alpha_n}", from=3-1, to=4-2]
\end{tikzcd}\right)
\end{equation}
Then, $P_{\mathbb{B}_{n+1}}'$ is a categorical idempotent
in $\Ch^-(\SBim_{\mathbb{B}_{n+1}})$ that is naturally homotopy equivalent to
$\mathbf{P}_{\mathbb{B}_{n+1}}$.
\end{proposition}
\noindent \textbf{Note:} The notation for the semi-infinite extension of
complexes in~(\ref{eq:nextstepconstruction}) is
explained in Section~\ref{ssec:constructpnplusone}.
\\ \\
While the construction of $P_{\mathbb{B}_{n+1}}'$ is in fact a special case of the
construction of categorified
projectors from categorical
diagonalization~\cite{elias_categorical_2017},
we will be presenting its construction in more elementary terms in this paper.
This recursive relation then immediately results in our theorem:
\begin{thm}
	\label{thm:homspaceparity}
	Given $n\geq 1$, the ring $H^\bullet\End(\mathbf{P}_{\mathbb{B}_n})$
	is supported
	in even cohomological degrees.
\end{thm}
\subsubsection{Proof by simultaneous induction}
The proof of Lemma~\ref{lem:eigenmaps} for $n=k$
depends on Theorem~\ref{thm:homspaceparity}
being true for $n=k-1$. Thus, we must prove all of these statements
simultaneously via induction. We will take $n=0$ as our base case, as
type $\mathbb{B}_1$ is identical to $\mathbb{A}_1$
and these results hold for type $\mathbb{A}$. We
will then proceed to assume that Theorem~\ref{thm:homspaceparity} holds for
$\mathbf{P}_{\mathbb{B}_{n}}$ and show that it implies the other statements for
type $\mathbb{B}_{n+1}$.
\subsubsection{Abbreviated notation}
In the following sections, we will denote the symmetrizer
complexes $P_{\mathbb{B}_n}$,
$\mathbf{P}_{\mathbb{B}_n}$, and $P_{\mathbb{B}_{n+1}}'$ by $P_n$,
$\mathbf{P}_{n}$, and $P_{n+1}'$ respectively. We
will also denote the generalized Jucys-Murphy complex
$\JM(\mathbb{B}_{n},\mathbb{B}_{n+1})$ by $\JM_{n+1}$, and a Soergel bimodule
$B_{s_{i_1}\dots s_{i_k}}$ corresponding to an element
$w=s_{i_1}\dots s_{i_k}$ of $W$ by $B_{i_1\dots i_k}$.
\subsection{From $\mathbb{A}_1=\mathbb{B}_1$ to $\mathbb{B}_2$}
As the simplified form of the complex $P_1\JM_2$ is slightly different
than when going from $\mathbb{B}_{2}$ to $\mathbb{B}_3$
and so on, we handle it separately here.
Recall from~\cite{elias_computation_2019} that Theorem~\ref{thm:homspaceparity}
holds for $n=1$, since type $\mathbb{B}_1$ is the same as $\mathbb{A}_1$.
We now compute
a simplified form for $P_1\JM_2$, which is
\[P_1F(s_2)F(s_1)F(s_2)F(s_2)F(s_1)F(s_2).\]
\subsubsection{Simplifying $P_1\JM_2$}
We will refer to sections~\ref{ssec:multtable}
and~\ref{ssec:subcomplexeliminations} to simplify our complexes. Since
$P_1B_1\simeq 0$, $P_1F(s_2)F(s_1)$ will be
\[\begin{tikzcd}[ampersand replacement=\&,cramped]
	{P_1B_2B_1} \& {t^{1/2}P_1B_2} \& tP_1
	\arrow[from=1-1, to=1-2]
	\arrow[from=1-2, to=1-3]
\end{tikzcd}\]
with the maps being the obvious counit maps. When we multiply by
$F(s_2)$ again, multiplication by $B_2$ will cause a term of
$P_1B_{212}$ to appear as in the complementary map to~(\ref{eqn:b2inclusion}),
while the remaining terms eliminate via~(\ref{eqn:3termelimination1}). The
complex after $F(s_2)F(s_1)F(s_2)$ has been tensored with will
be
\[\begin{tikzcd}[ampersand replacement=\&,cramped]
	\&\& {P_1B_{212}} \\
	{t^{1/2}P_1B_2B_1} \& {tP_1B_2} \& {t^{3/2}P_1}
	\arrow[from=1-3, to=2-1]
	\arrow[from=2-1, to=2-2]
	\arrow[from=2-2, to=2-3]
\end{tikzcd}\]
After tensoring with $F(s_2)$ again, the complex
$P_1F(s_2)F(s_1)F(s_2)F(s_2)$
will become
\[\begin{tikzcd}[ampersand replacement=\&,cramped]
	\&\&\& {q^{1/2}P_1B_{212}} \\
	\\
	{t^{1/2}P_1B_{212}} \& {tP_1B_2B_1} \& {t^{3/2}P_1B_2} \& {t^{2}P_1}
	\arrow["\xi_1"{description}, from=1-4, to=3-1]
	\arrow[from=3-1, to=3-2]
	\arrow[from=3-2, to=3-3]
	\arrow[from=3-3, to=3-4]
\end{tikzcd}\]
Where the map $\xi_1$ is
\[\q[scale=0.5]{
		\qbox[n=4,s=4,N={},S={},hlen=4,id=idempbox-top] \n
		\qspace \qspace \qspace \qspace
		\qbox[n=1,s=1,N={},S={},name=$\alpha_2$,frame color=none,hlen=1]\n
		\qbox[n=4,s=4,N={},S={},hlen=4,id=idempbox-bottom]
		\qwire[dom=n-1-idempbox-bottom,cod=s-1-idempbox-top,color=orange]
		\qwire[dom=n-2-idempbox-bottom,cod=s-2-idempbox-top,color=blue]
		\qwire[dom=n-4-idempbox-bottom,cod=s-4-idempbox-top,color=orange]
		\qwire[dom=n-1-idempbox-top,color=orange]
		\qwire[dom=n-2-idempbox-top,color=blue]
		\qwire[dom=n-4-idempbox-top,color=orange]
		\qwire[cod=s-1-idempbox-bottom,color=orange]
		\qwire[cod=s-2-idempbox-bottom,color=blue]
		\qwire[cod=s-4-idempbox-bottom,color=orange]
}
-
\q[scale=0.5]{
		\qbox[n=4,s=4,N={},S={},hlen=4,id=idempbox-top] \n
		\qspace \qspace
		\qbox[n=1,s=1,N={},S={},name=$\alpha_2$,frame color=none,hlen=1]\n
		\qbox[n=4,s=4,N={},S={},hlen=4,id=idempbox-bottom]
		\qwire[dom=n-1-idempbox-bottom,cod=s-1-idempbox-top,color=orange]
		\qwire[dom=n-2-idempbox-bottom,cod=s-2-idempbox-top,color=blue]
		\qwire[dom=n-4-idempbox-bottom,cod=s-4-idempbox-top,color=orange]
		\qwire[dom=n-1-idempbox-top,color=orange]
		\qwire[dom=n-2-idempbox-top,color=blue]
		\qwire[dom=n-4-idempbox-top,color=orange]
		\qwire[cod=s-1-idempbox-bottom,color=orange]
		\qwire[cod=s-2-idempbox-bottom,color=blue]
		\qwire[cod=s-4-idempbox-bottom,color=orange]
}.\]
where the blue strands correspond to $B_1$ and the orange strands to
$B_2$, and the boxes are the projection/inclusion for $B_{212}$. The simplifications
involved are the same as in the $n > 1$ case and will be handled in greater detail
in the proof for Lemma~\ref{lem:pjmsimplified}.
\\ \\
Now, when tensoring
with $F(s_1)$, $B_1$ will eliminate the entire bottom row
via~(\ref{eqn:3termelimination2}) and take the top term to
$P_1B_{2}B_1$ via~(\ref{eqn:b21downstepiso}), so the complex
$P_1F(s_2)F(s_1)F(s_2)F(s_2)F(s_1)$ will be homotopy
equivalent to
\[\begin{tikzcd}[ampersand replacement=\&,cramped]
	\&\& {q^{1/2}P_1B_2B_1} \& {q^{1/2}t^{1/2}P_1B_{212}} \\
	\\
	{tP_1B_{212}} \& {t^{3/2}P_1B_2B_1} \& {t^2P_1B_2} \& {t^{5/2}P_1}
	\arrow[from=1-3, to=1-4]
	\arrow["\xi_1"{description}, from=1-4, to=3-1]
	\arrow[from=3-1, to=3-2]
	\arrow[from=3-2, to=3-3]
	\arrow[from=3-3, to=3-4]
\end{tikzcd}\]
When tensoring with the $B_2$ from $F(s_2)$, all terms except the middle
$q^{1/2}t^{1/2}P_1B_{212}$ term
(and the lowest-$t$-degree $q^{1/2}P_1B_2B_1$, which
will be taken to $q^{1/2}PB_2$) will eliminate via the relations
in sections~\ref{ssec:multtable}
and~\ref{ssec:subcomplexeliminations}, with the remainder eliminating
via the observation~(\ref{eqn:finaleliminationmap}) in the proof of
Lemma~\ref{lem:pjmsimplified}. So, $P_1\JM_2$ is homotopy equivalent to
\begin{equation}
	\label{eqn:p1jm2simplified}
	\begin{tikzcd}[ampersand replacement=\&,cramped]
	\& {q^{1/2}P_1B_2} \& {q^{1/2}t^{1/2}P_1B_2B_1} \& {q^{1/2}tP_1B_{212}} \\
	\\
	{t^{3/2}P_1B_{212}} \& {t^{2}P_1B_2B_1} \& {t^{5/2}P_1B_2} \& {t^{3}P_1}
	\arrow[from=1-2, to=1-3]
	\arrow[from=1-3, to=1-4]
	\arrow["\xi_1"{description}, from=1-4, to=3-1]
	\arrow[from=3-1, to=3-2]
	\arrow[from=3-2, to=3-3]
	\arrow[from=3-3, to=3-4]
\end{tikzcd}\end{equation}
\subsubsection{Existence of the eigenmaps}
Since $P_1$ is a unital idempotent, there is a homotopy equivalence
$\Hom(P_1,P_1\JM_2)\simeq \Hom(\mathbb{1},P_1\JM_2)$. So, we apply the
$\Hom(\mathbb{1},-)$ functor to~(\ref{eqn:p1jm2simplified}) and, after observing that
\[\Hom(\mathbb{1},P_1B_2B_1)\simeq \Hom(B_1,P_1B_2)\simeq 0\]
via the properties of rigid monoidal categories and
Theorem~\ref{thm:idempotentmorphismzero}, and
\[\Hom(\mathbb{1},P_1B_{212})\simeq 0\]
by Lemma~\ref{lem:homspacezerolightleaves},
we are left with
\begin{align*}
	\Hom &(\mathbb{1},q^{1/2}P_1B_{2}) \\
	     &{\big\downarrow}\\
	     \biggl( \Hom(\mathbb{1},t^{5/2}P_1B_{2})
	     & \xrightarrow{\Hom(\mathbb{1},\id_{P_1} \otimes
	     \q[scale=0.4]{\qcounit[color=black]})}
     \Hom(\mathbb{1},t^{3}P_1) \biggr)
\end{align*}
Bt the double leaves basis, any morphism $\mathbb{1}\to P_1B_2$ is
a morphism $\mathbb{1}\to P_1$ tensored with a unit for $B_2$, so we have
natural isomorphisms
\[\Hom(\mathbb{1},P_1B_2)\simeq q^{1/2}\Hom(\mathbb{1},P_1)\simeq
\Hom(\mathbb{1},P_1)\otimes \Hom(\mathbb{1},B_2).\]
Therefore, the arrow
on the bottom row is an injection (multiplying by $\alpha_2$). Combining
this fact and the fact that $\End(P_1)$ is already known to be parity
in homology, the proof
of Lemma~\ref{lem:homspacesurjection} can be followed to yield a map
$\beta_1:qP_1\to P_1\JM_2$ of the desired form, after which
Theorem~\ref{thm:homspaceparity} can be proven for $n=2$ by simply
following the steps in section~\ref{ssec:constructpnplusone} and onwards.
\subsection{From $\mathbb{B}_n$ to $\mathbb{B}_{n+1}$ for $n > 1$}
The proof of
Lemma~\ref{lem:eigenmaps}
for higher values of $n$ is much like the previous case, except for
the presence of more intermediate terms in the simplified form of
$P_n\JM_{n+1}$.
\subsubsection{The
$P_nB_{n+1}$-multiplication table}
\label{ssec:multtable}
In this section, we supply calculations which will help us simplify the complex
$P_n\JM_{n+1}$. Throughout, we will use the shorthands
\[P_{n}B_{n+1}\dots B_i:= P_n\otimes \bigotimes_{k=1}^{n+1-i}B_{n+1-k}\]
and
\[P_nB_{n+1}\dots B_{212}\dots B_{i}:= P_nB_{n+1}\dots B_3 B_{212}
\otimes \bigotimes_{k=3}^i B_k\]
where $n+1 \geq i \geq 3$.
\begin{lemma}
	\label{lem:basicmultiplications1}
Given an elementary Bott-Samelson $B_j$, we have the following
	homotopy equivalences:
	\begin{enumerate}[label=(\arabic*)]
		\item If $|i-j| >1$, then
			\begin{equation}
				\label{eqn:zerorelation1}
				P_nB_{n+1}\dots B_iB_j\simeq
				0\simeq P_nB_{n+1}\dots B_{212}\dots
				B_iB_j
			\end{equation}
		\item If $j=i+1$, then
			\begin{equation*}
			P_nB_{n+1}\dots B_{212}\dots B_iB_j=
			P_nB_{n+1}\dots B_{212}\dots B_{i+1}
			\end{equation*}
			and
			\begin{equation*}
			P_nB_{n+1}\dots B_iB_j\simeq
			P_nB_{n+1}\dots B_{i+1}
			\end{equation*}
			\item If $j=i-1$, then
		\begin{equation*}
		P_nB_{n+1}\dots B_iB_j= P_nB_{n+1}\dots B_{i-1}
		\end{equation*}
		when $i \geq 2$ and 
		\begin{equation*}
		P_nB_{n+1}\dots B_{212}\dots B_iB_j\simeq
		P_nB_{n+1}\dots B_{212}\dots B_{i-1}
		\end{equation*}
		when $i > 3$ with the homotopy equivalence being the same
		as in (\ref{eqn:downshiftiso}) or
		\begin{equation*}
		P_nB_{n+1}\dots B_{212}B_3B_2 \simeq P_nB_{n+1}\dots
		B_{212}
		\end{equation*}
		when $i=3$.
	\item If $j=i$, then
	\begin{align*}
		\begin{split}
		P_nB_{n+1}\dots B_iB_j &\simeq (q^{1/2}\mathbb{1}
		\oplus q^{-1/2}\mathbb{1})P_nB_{n+1}\dots B_i\\
		P_nB_{n+1}\dots B_{212}\dots B_iB_j &\simeq
		(q^{1/2}\mathbb{1} \oplus q^{-1/2}\mathbb{1})
		P_nB_{n+1}\dots B_{212}\dots B_i
	\end{split}
	\end{align*}
	where the inclusion of each direct
	summand is given by $\id_{P_nB_{n+1}\dots
	B_{i+1}}\otimes \q[scale=0.2]{
	\qwire[color=purple] \qbox[hlen=1,n=0,s=0] \qwire[color=purple]\n
	\qcupraw[hlen=3, frame color=purple]\n
	\qspace \qwire[vlen=1.4,color=purple]}$
	and 
	$\id_{P_nB_{n+1}\dots
	B_{212}\dots B_{i-1}}\otimes \q[scale=0.2]{
	\qwire[color=purple] \qbox[hlen=1,n=0,s=0] \qwire[color=purple]\n
	\qcupraw[hlen=3, frame color=purple]\n
	\qspace \qwire[vlen=1.4,color=purple]}$ in both cases,
	with the value in the box being
	$1$ for the inclusion of the $q^{-1/2}$-shifted term and
	$\alpha_i/2$ for the inclusion of the $q^{1/2}$-shifted term. The diagram is
	flipped and the values in the box are swapped for the projections.
	\end{enumerate}
\end{lemma}
\begin{proof}
We first focus on the case of
$P_nB_{n+1}\dots B_iB_j$.
\begin{itemize}
	\item If $j < i-1$, then
\[P_nB_{n+1}\dots B_iB_j\simeq P_nB_jB_{n+1}\dots B_i\simeq 0\]
which proves one part of (1).
\item If $j=i-1$, (3) obviously holds.
\item If $j=i$, relation (4) obviously holds by the quadratic relation in
$\SBim$.
\item If $j= i+1$, then
\[P_nB_{n+1}\dots B_{i+1} B_i B_{i+1}\simeq
P_nB_{n+1}\dots B_{i+2} (B_{(i+1)i(i+1)}\oplus B_{i+1}).\]
However, since $i \geq 3$, $B_{(i+1)i(i+1)}=B_{i(i+1)i}$ is a direct summand
of $B_{i}B_{i+1}B_i$, meaning the first direct summand is
\[P_nB_{n+1}\dots B_{i+2}B_{i}B_{i+1}B_i\simeq 0\]
which has already been proven to be contractible. Therefore, (2) holds,
and the homotopy equivalence is simply the inclusion of the direct summand
given as follows:
			\begin{equation}
	\label{eqn:downshiftiso}
	P_nB_{n+1}\dots \textcolor{purple}
	{B_{i+1}} \xrightarrow{\id_{P_{n}B_{n+1}\dots
		B_{i+2}}\otimes \q[scale=0.15]{\qwire[color=purple]
		\qunit[color=pink,
		frame color=pink, id=u] \qwire[color = purple]
		\n
					\qcupraw[hlen=3, frame color=purple]\n
			\qspace \qwire[vlen=1.4, color=purple]
		\qwire[dom=C-u, cod=N-u, color=pink]
		}
	} P_nB_{n+1}\dots \textcolor{purple}{B_{i+1}}\textcolor{pink}{B_{i}}
	\textcolor{purple}{B_{i+1}}.
	\end{equation}
\item If $j > i+1$, then
	\begin{align*}P_nB_{n+1}\dots B_iB_j &\simeq
	P_nB_{n+1}\dots B_{j+1}B_jB_{j-1}B_jB_{j-2}\dots B_i
	\\ &\simeq
P_nB_{n+1}\dots B_{j+1}(B_j\oplus B_{j(j-1)j} )B_{j-2}\dots B_i \simeq 0
\end{align*}
where the first direct summand is eliminated by $B_{j-2}$  (which is at least
$B_i$) from the right and the second direct summand is contractible
by the relation given in (2). This proves the other part of (1).
\item For the special case of (3) where
	$i=3$, we see that $B_{212}B_3B_2\simeq 
	B_{212}\oplus B_{21232}$ by Hecke category diagrammatics with
	the
	inclusion of $B_{212}$ given by the map in~(\ref{eqn:b2123downstepiso}). By
	the relations in the type $B$ coxeter group,
	$B_{21232}=B_{21323}=B_{23123}$, which is a direct summand of
	$B_2B_3B_1B_2B_3$. Thus,
	\begin{align*}P_nB_{n+1}\dots B_3B_{21232} &\hookrightarrow
	P_nB_{n+1}\dots B_{3}B_2B_3B_1B_2B_3 \\
	&\simeq P_nB_{n+1}\dots B_3B_1B_2B_3 \simeq 0
\end{align*}
	with the contraction happening due to the part of (1) previously proven. Thus,
	the remaining direct summand is $P_nB_{n+1}\dots B_{212}$,
	and the homotopy equivalence is given by
\begin{equation}
	\label{eqn:b2123downstepiso}
P_nB_{n+1}\dots B_{\textcolor{orange}{2}\textcolor{blue}{1}\textcolor{orange}{2}}
	\xrightarrow{
	\id_{P_nB_{n+1}\dots B_3}\otimes 
	\q[scale=0.25]{
	\qbox[n=3,s=3,N={},S={},hlen=3,id=idempbox2] 
\qunit[color=green, frame color=green,
id=newu]
\qwire[color=orange, id=splitwireend]
		\n
		\qwire[color=orange] \qwire[color=blue]\qwire[color=orange,
		id=splitwire]
		\n
		\qbox[n=3,s=3,N={},S={},hlen=3,id=idempbox] 
\qwire[dom=n-3-idempbox,color=orange,vlen=0.25]
\qwire[dom=n-2-idempbox,color=blue,vlen=0.25]
\qwire[dom=n-1-idempbox,color=orange,vlen=0.25]
\qwire[cod=s-1-idempbox,color=orange,vlen=0.25]
\qwire[cod=s-3-idempbox,color=orange,vlen=0.25]
\qwire[cod=s-2-idempbox,color=blue,vlen=0.25]
\qwire[cod=s-1-idempbox,color=orange,vlen=0.25]
\qwire[dom=n-3-idempbox2,color=orange,vlen=0.25]
\qwire[dom=n-2-idempbox2,color=blue,vlen=0.25]
\qwire[dom=n-1-idempbox2,color=orange,vlen=0.25]
\qwire[cod=s-1-idempbox2,color=orange,vlen=0.25]
\qwire[cod=s-3-idempbox2,color=orange,vlen=0.25]
\qwire[cod=s-2-idempbox2,color=blue,vlen=0.25] 
\qwire[cod=s-1-idempbox2,color=orange,vlen=0.25]
\qwire[dom=C-newu, cod=N-newu, color=green] 
	\qwire[dom=C-splitwire, cod=S-splitwireend, color=orange]
}
}
	P_nB_{n+1}P_nB_{n+1}\dots B_{\textcolor{orange}{2}
	\textcolor{blue}{1}\textcolor{orange}{2}}\textcolor{green}{B_3}
	\textcolor{orange}{B_2}
\end{equation}
where the boxes are the projection/inclusion for $B_{212}$.
\end{itemize}
We now move on to complexes of the form
$P_nB_{n+1}\dots B_{212}\dots B_iB_j$, where $i \geq 3$.
\begin{itemize}
	\item If $j > i+1$, then
		\[P_nB_{n+1}\dots B_{212}\dots B_iB_j\simeq
		P_nB_{n+1}\dots B_3B_jB_{212}\dots B_i \simeq 0\]
		since $j > 4$ by the part of (1) previously proven.
	\item If $j = i+1$, the complex is not reduced or contracted
		in any way.
	\item If $j=i$, the quadratic relation proves (4) as usual.
	\item If $j = i-1$, we assume $i > 3$ (as we have already proven the
		case where $i =3$. Then,
		\begin{align*}
		P_nB_{n+1}\dots B_{212}\dots B_iB_j &\simeq
		P_nB_{n+1}\dots B_{212} \dots B_{i-1}B_iB_{i-1}\\
		&\simeq P_{n}B_{n+1}\dots B_{212}\dots (B_{i-1}\oplus
		B_{(i-1)i(i-1)}).
\end{align*}
Now, since $B_{(i-1)i(i-1)}=B_{i(i-1)i}$ is a direct summand of
$B_iB_{i-1}B_i$, the second direct summand will be a direct summand of
\[P_nB_{n+1}\dots B_{212}\dots B_iB_{i-1}B_i\simeq 0\]
where the contraction occurs as in the case where $j > i+1$. Thus, the surviving
term is $P_nB_{n+1}\dots B_{212}\dots B_{i-1}$ and the relation in (3) holds.
\item If $j < i-1$ but $j >2$, then 
	\begin{align*}
	P_nB_{n+1}\dots B_{212}\dots B_iB_j &\simeq 
	P_nB_{n+1}\dots B_{212}\dots B_jB_{j+1}B_j \dots B_i \\
					    &\simeq P_nB_{n+1}\dots
					    B_{212}\dots (B_j \oplus
					    B_{j(j+1)j})\dots B_i
	\end{align*}
	where both direct summands are contractible by the case of (1)
	proven before. If $j=2$, then the special case of (3) gives that
	\begin{align*}
		P_nB_{n+1}\dots B_{212}\dots B_iB_2 &\simeq
		P_nB_{n+1}\dots B_{212}B_3B_2\dots B_i \\
		&\simeq P_nB_{n+1}\dots B_{212}B_4 \dots B_i
	\end{align*}
	which will be contractible by the previous results.
	\\ \\
	If
	$j=1$, then relation (6) in
	Lemma~\ref{lem:basicmultiplications2} (which is proven independently)
	will make the complex homotopy equivalent to
	\[P_nB_{n+1}\dots B_{2}B_1 B_3 \dots B_i\]
	which is again contractible.
\end{itemize}
\end{proof}
\noindent Additionally, there are the following identities involving
$P_nB_{n+1}\dots B_{212}$ and $P_nB_{n+1}\dots B_{2}B_1$.
\begin{lemma}
	\label{lem:basicmultiplications2}
	The following identities hold:
\begin{enumerate}[label=(\arabic*)]
	\setcounter{enumi}{4}
	\item Given the complex $P_nB_{n+1}\dots B_{2}B_1$, tensoring with
	$B_2$ on the right yields
\[P_nB_{n+1}\dots B_{2}B_1B_2\simeq
P_nB_{n+1}\dots B_{212} \oplus P_nB_{n+1}\dots B_{2}\]
\item Given the complex $P_nB_{n+1}\dots B_{212}$, tensoring with $B_1$ on the
	right yields
	\[P_nB_{n+1}\dots B_{212}B_1\simeq
	P_nB_{n+1}\dots B_2 B_1\]
	\end{enumerate}
\end{lemma}
\begin{proof}
We prove these relations in order.
	\begin{itemize}
	\item Relation (5) holds by the fact that in type $B$,
$B_{212}\neq B_{121}$ as an indecomposable Soergel bimodule and is thus not
contracted by the remainder of the complex, and the morphisms
affording the homotopy equivalence are simply the
inclusions of the direct summands for $B_2B_1B_2 \simeq B_2\oplus B_{212}$.
The inclusion of $P_nB_{n+1}\dots B_2$ is given by
\begin{equation}
	\label{eqn:b2inclusion}
	P_nB_{n+1}\dots \textcolor{orange}
	{B_{2}} \xrightarrow{\id_{P_{n}B_{n+1}\dots
		B_{i+2}}\otimes \q[scale=0.15]{\qwire[color=orange]
		\qunit[color=blue,
		frame color=blue, id=u] \qwire[color = orange]
		\n
					\qcupraw[hlen=3, frame color=orange]\n
			\qspace \qwire[vlen=1.4, color=orange]
		\qwire[dom=C-u, cod=N-u, color=blue]
		}
	} P_nB_{n+1}\dots \textcolor{orange}{B_{2}}\textcolor{blue}{B_{1}}
	\textcolor{orange}{B_{2}}.
	\end{equation}
	and that of $P_n B_{n+1}\dots B_{212}$ is given by tensoring in the projector
	for $B_{212}$ instead.
\item For (6), we observe that in type
	$B$, $B_{212}B_1\simeq B_{2121}\oplus B_{21}
	=B_{1212}\oplus B_2B_1$ from
	Hecke category diagrammatics, with
	the inclusion of $B_{21}$ being given by the morphism
\begin{equation}
	\label{eqn:b21downstepiso}
	P_nB_{n+1}\dots \textcolor{orange}{B_{2}}\textcolor{blue}{B_1}
	\xrightarrow{\id_{P_nB_{n+1}\dots B_3} \otimes
		\q[scale=0.15]{
		\qbox[n=3,s=3,N={},S={},hlen=3,id=idempbox] 
		\qwire[color=blue] \n
		\qwire[color=orange]
		\qwire[color=blue]
		\qunit[color=orange,frame color=orange, id=u]
		\qwire[color=blue]
		\n
		\qwire[color=orange] \qcupraw[frame color=blue, hlen=3]\n
		\qwire[color=orange] \qspace \qwire[color=blue,vlen=1.4]
\qwire[dom=n-1-idempbox,color=orange,vlen=0.25]
\qwire[dom=n-3-idempbox,color=orange,vlen=0.25]
\qwire[dom=n-2-idempbox,color=blue,vlen=0.25]
\qwire[cod=s-1-idempbox,color=orange,vlen=0.25]
\qwire[cod=s-3-idempbox,color=orange,vlen=0.25]
\qwire[cod=s-2-idempbox,color=blue,vlen=0.25]
\qwire[dom={C-u}, cod={S-3-idempbox}, color=orange]
		}
	}
	P_nB_{n+1}\dots B_{\textcolor{orange}{2}
	\textcolor{blue}{1}\textcolor{orange}{2}}\textcolor{blue}{B_1}
\end{equation}
where the box is the projector for $B_{212}$.
	Then,
	\[P_nB_{n+1}\dots B_{212}B_1\simeq
	P_nB_{n+1}\dots B_3(B_{2121}\oplus B_{2}B_1)\simeq P_nB_{n+1}\dots B_2B_1\]
	as the first direct summand $P_nB_{n+1}\dots B_3B_{2121}$ is a direct summand
	of $P_nB_{n+1}\dots B_{3}B_1B_2B_1B_2$, which is zero by
	relation (1) in Lemma~\ref{lem:basicmultiplications1}. Thus, the
	morphism in~(\ref{eqn:b21downstepiso}) is a homotopy equivalence.
\end{itemize}
\end{proof}
\subsubsection{Subcomplex eliminations}
\label{ssec:subcomplexeliminations}
We use the homotopy equivalences in
Lemmas~\ref{lem:basicmultiplications1}~and~\ref{lem:basicmultiplications2} to
show that certain complexes are eliminated by tensoring with certain Bott-Samelsons
on the right.
\begin{lemma}
	\label{lem:subcomplexeliminations}
	The following homotopy equivalences hold:
	\begin{enumerate}[label=(\arabic*)]
		\item Suppose a complex $C$ of the form\\\\
{\centering
$
{\begin{aligned}
	P_nB_{n+1} B_n \dots &\textcolor{pink}{B_{i+1}} B_{i}
		\textcolor{purple}{B_{i-1}}
	\\
			     &{\big\downarrow}
			     \q[scale=0.5]{\qwire[color=pink] \qwire
			     \qcounit[color=purple, frame color=purple, id=u]
		     \qwire[dom=C-u, cod=S-u, color=purple]}
	\\
	t^{1/2}P_nB_{n+1}B_n\dots &\textcolor{pink}{B_{i+1}} B_i
	\\
&{\big\downarrow} \q[scale=0.5]{\qwire[color=pink]
\qcounit[color=black]}
	\\
	t P_n B_{n+1}B_n\dots &\textcolor{pink}{B_{i+1}}
\end{aligned}}$
\hspace{0.05\textwidth}or\hspace{0.05\textwidth}
${\begin{aligned}
	P_nB_{n+1} B_n \dots B_{212}\dots  &\textcolor{pink}{B_{i-1}} B_{i}
		\textcolor{purple}{B_{i+1}}
	\\
			     &{\big\downarrow}
			     \q[scale=0.5]{\qwire[color=pink] \qwire
			     \qcounit[color=purple, frame color=purple, id=u]
		     \qwire[dom=C-u, cod=S-u, color=purple]}
	\\
	t^{1/2}P_nB_{n+1}B_n\dots B_{212}\dots&\textcolor{pink}{B_{i-1}} B_i
	\\
&{\big\downarrow} \q[scale=0.5]{\qwire[color=pink]
\qcounit[color=black]}
	\\
	t P_n B_{n+1}B_n\dots B_{212}\dots &\textcolor{pink}{B_{i-1}}
\end{aligned}}$
}\\\\
Then, $CB_i\simeq 0$.
\item Suppose a complex $C$ of the form
\begin{align*}
	P_nB_{n+1}\dots	&B_{\textcolor{orange}{2}\textcolor{blue}{1}\textcolor{orange}{2}}\\
&{\big\downarrow} \q[scale=0.25]{\qwire[color=orange]\qwire[color=blue]
			     \qcounit[color=orange, frame color=orange, id=u]
			     \n\qbox[hlen=3,N={},S={},n=3,s=3,
			     id=idempbox]
		     \qwire[dom=C-u, cod=S-u, color=orange]
\qwire[dom=n-1-idempbox,color=orange,vlen=0.25]
\qwire[dom=n-3-idempbox,color=orange,vlen=0.25]
\qwire[dom=n-2-idempbox,color=blue,vlen=0.25]
\qwire[cod=s-1-idempbox,color=orange,vlen=0.25]
\qwire[cod=s-3-idempbox,color=orange,vlen=0.25]
\qwire[cod=s-2-idempbox,color=blue,vlen=0.25]}	  
	\\
	t^{1/2}P_nB_{n+1}\dots &\textcolor{orange}{B_2}\textcolor{blue}{B_1}\\
			       &{\big\downarrow}
			       \q[scale=0.25]{\qwire[color=orange]
				       \qcounit[color=blue,frame color=blue,id=u]
			       \qwire[dom=C-u, cod=S-u, color=blue]}
	\\
	tP_nB_{n+1}\dots &\textcolor{orange}{B_{2}}.
\end{align*}
Then, $CB_1\simeq 0$.
\item Suppose a complex $C$ of the form
	\begin{align*}
	\begin{split}
	\label{eqn:5termelimination}
	P_nB_{n+1}\dots &\textcolor{green}{B_3}
B_{\textcolor{orange}{2}\textcolor{blue}{1}\textcolor{orange}{2}}
\textcolor{green}{B_3}
\\
&{\big\downarrow} \q[scale=0.25]{
	\qwire[color=green]\qbox[hlen=3,N={},S={},n=3,s=3,id=idempbox]
			    \qcounit[color=green, frame color=green, id=u]
		     \qwire[dom=C-u, cod=S-u, color=green]
\qwire[dom=n-1-idempbox,color=orange,vlen=0.25]
\qwire[dom=n-3-idempbox,color=orange,vlen=0.25]
\qwire[dom=n-2-idempbox,color=blue,vlen=0.25]
\qwire[cod=s-1-idempbox,color=orange,vlen=0.25]
\qwire[cod=s-3-idempbox,color=orange,vlen=0.25]
\qwire[cod=s-2-idempbox,color=blue,vlen=0.25]}	 
\\
	t^{1/2}P_nB_{n+1}\dots &\textcolor{green}{B_3}
B_{\textcolor{orange}{2}\textcolor{blue}{1}\textcolor{orange}{2}}
\\
		&{\big\downarrow} \q[scale=0.25]{\qwire[color=green]
	\qwire[color=orange]\qwire[color=blue]
			     \qcounit[color=orange, frame color=orange, id=u]
			     \n \qwire[color=green]
			     \qbox[hlen=3,N={},S={},n=3,s=3,
			     id=idempbox]
		     \qwire[dom=C-u, cod=S-u, color=orange]
\qwire[dom=n-1-idempbox,color=orange,vlen=0.25]
\qwire[dom=n-3-idempbox,color=orange,vlen=0.25]
\qwire[dom=n-2-idempbox,color=blue,vlen=0.25]
\qwire[cod=s-1-idempbox,color=orange,vlen=0.25]
\qwire[cod=s-3-idempbox,color=orange,vlen=0.25]
\qwire[cod=s-2-idempbox,color=blue,vlen=0.25]}	 
\\
tP_nB_{n+1}\dots &\textcolor{green}{B_3}\textcolor{orange}{B_2}
\textcolor{blue}{B_1}\\
&{\big\downarrow} \q[scale=0.3]{\qwire[color=green]
	\qwire[color=orange]
	\qcounit[color=blue, frame color=blue,id=u]
			 		     \qwire[dom=C-u, cod=S-u, color=blue]}	
\\
t^{3/2}P_nB_{n+1}\dots & \textcolor{green}{B_3}\textcolor{orange}{B_2} \\
&{\big\downarrow} \q[scale=0.3]{\qwire[color=green]
	\qcounit[color=orange, frame color=orange,id=u]
\qwire[dom=C-u, cod=S-u, color=orange]}	
\\
t^2P_nB_{n+1}\dots &\textcolor{green}{B_3}.
\end{split}
\end{align*}
Then, $CB_2\simeq 0$.
\end{enumerate}
Versions of these complexes where the directions of the arrows are reversed
(with appropriate changes to grading shifts) and the diagrams are flipped satisfy
the same relations.
\end{lemma}
\begin{proof}
	We prove the contractions in sequence.
	\begin{enumerate}[label=(\arabic*)]
		\item We demonstrate the proof only on the first
			of the two complexes, since all relevant isomorphisms
			only involve the rightmost three terms and are thus
			identical. When tensored on
			the right with $B_i$, it yields the following complex
after the quadratic relation in $\SBim$ and~(\ref{eqn:downshiftiso}):
\begin{equation}
	\label{eqn:3termelimination1}
\begin{tikzcd}[ampersand replacement=\&,cramped]
	{P_nB_{n+1}\dots \textcolor{pink}{B_{i+1}}B_i} \&\& \\
	\q[scale=0.3]{ \qwire[color=pink] \qspace \qwire[id=outwire] \n
		\qwire[color=pink] \qcapraw[hlen=3, id=capup] \n
		\qwire[color=pink] \qwire
		\qcounit[color=purple, frame color=purple, id=cou] \qwire \n
		\qwire[color=pink] \qwire \qwire[color=purple]\qwire \n
		\qwire[color=pink]\qwire
		\qunit[color=purple, frame color=purple,
		id=u]
		\qwire
		
		\n
		\qwire[color=pink] \qcupraw[hlen=3, id=cupdown] \n
		\qwire[color=pink] \qspace \qwire[id=inwire]
		\qwire[dom=C-cou, cod=S-cou, color=purple]
		\qwire[dom=C-u, cod=N-u, color=purple]
		\qwire[cod=S-outwire, vlen=0.4]
		\qwire[dom=N-inwire,vlen=0.4]
} \\
	{t^{1/2}q^{1/2}P_nB_{n+1}\dots \textcolor{pink}{B_{i+1}}B_i} \&
	\bigoplus \& {t^{1/2}q^{-1/2}P_nB_{n+1}\dots \textcolor{pink}{B_{i+1}}B_i} \\
	\&\& \q[scale=0.3]{ \qwire[color=pink]  \qspace \qspace \qwire[id=outwire]\n
		\qwire[color=pink]\qcounit[color=black]\qspace
		\qwire	
		\n
		\qwire[color=pink] \qcupraw[hlen=3, id=cupdown] \n
		\qwire[color=pink] \qspace \qwire[id=inwire]
		\qwire[dom=N-inwire,vlen=0.4]
} \\
	\&\& {tP_nB_{n+1}\dots \textcolor{pink}{B_{i+1}}B_i}
	\arrow[shift right=10, from=1-1, to=3-1]
	\arrow[shift right=10, from=3-3, to=5-3]
\end{tikzcd}
\end{equation}
where both of the arrows pictured will be identities, making the entire complex
contractible and eliminating all three terms.
\item Tensor the complex on the right with $B_1$. After applying
equation~\ref{eqn:b21downstepiso}, the complex becomes of the form
\begin{equation}
\label{eqn:3termelimination2}
\begin{tikzcd}[ampersand replacement=\&,cramped]
	{P_nB_{n+1}\dots \textcolor{orange}{B_2}\textcolor{blue}{B_1}} \&\& \\
	\q[scale=0.3]{
		\qwire[color=orange] \qspace \qwire[color=blue,id=outwire]\n
		\qwire[color=orange] \qcapraw[frame color=blue,hlen=3]\n
		\qwire[color=orange] \qwire[color=blue]\qcounit[color=orange,
		frame color=orange, id=cou] \qwire[color=blue]
		\n
	\qbox[n=3,s=3,N={},S={},hlen=3,id=idempbox] 
		\qwire[color=blue] \n
		\qwire[color=orange]
		\qwire[color=blue]
		\qunit[color=orange,frame color=orange, id=u]
		\qwire[color=blue]
		\n
		\qwire[color=orange] \qcupraw[frame color=blue, hlen=3]\n
		\qwire[color=orange] \qspace \qwire[color=blue,vlen=1.4]
\qwire[dom=n-1-idempbox,color=orange,vlen=0.25]
\qwire[dom=n-3-idempbox,color=orange,vlen=0.25]
\qwire[dom=n-2-idempbox,color=blue,vlen=0.25]
\qwire[cod=s-1-idempbox,color=orange,vlen=0.25]
\qwire[cod=s-3-idempbox,color=orange,vlen=0.25]
\qwire[cod=s-2-idempbox,color=blue,vlen=0.25]
\qwire[cod=S-outwire,vlen=0.4,color=blue]
\qwire[dom={C-u}, cod={S-3-idempbox}, color=orange]`
\qwire[cod={C-cou}, dom={N-3-idempbox}, color=orange]

	}
	\\
	{P_nB_{n+1}\dots \textcolor{orange}{B_2}\textcolor{blue}{B_1}} \& \bigoplus \& {P_nB_{n+1}\dots \textcolor{orange}{B_2}\textcolor{blue}{B_1}} \\
	\&\& 
	\q[scale=0.3]{ \qwire[color=orange]  \qspace \qspace \qwire[color=blue,
		id=outwire]\n
		\qwire[color=pink]\qcounit[color=blue,frame color=blue,
		id=u]\qspace
		\qwire	
		\n
		\qwire[color=orange] \qcupraw[hlen=3, frame color=blue,id=cupdown] \n
		\qwire[color=orange] \qspace \qwire[color=blue, id=inwire]
		\qwire[dom=N-inwire, color=blue,vlen=0.4]
		\qwire[dom=C-u, cod=S-u, color=blue]
}
	\\
	\&\& {P_nB_{n+1}\dots \textcolor{orange}{B_2}\textcolor{blue}{B_1}}
	\arrow[shift right=10, from=1-1, to=3-1]
	\arrow[shift right=10, from=3-3, to=5-3]
\end{tikzcd}
\end{equation}
where the rightmost arrow is obviously an isomorphism, and the other arrow
is an isomorphism because
\[
	\q[scale=0.3]{
		\qwire[color=orange] \qspace \qwire[color=blue,id=outwire]\n
		\qwire[color=orange] \qcapraw[frame color=blue,hlen=3]\n
		\qwire[color=orange] \qwire[color=blue]\qcounit[color=orange,
		frame color=orange, id=cou] \qwire[color=blue]
		\n
	\qbox[n=3,s=3,N={},S={},hlen=3,id=idempbox] 
		\qwire[color=blue] \n
		\qwire[color=orange]
		\qwire[color=blue]
		\qunit[color=orange,frame color=orange, id=u]
		\qwire[color=blue]
		\n
		\qwire[color=orange] \qcupraw[frame color=blue, hlen=3]\n
		\qwire[color=orange] \qspace \qwire[color=blue,vlen=1.4]
\qwire[dom=n-1-idempbox,color=orange,vlen=0.25]
\qwire[dom=n-3-idempbox,color=orange,vlen=0.25]
\qwire[dom=n-2-idempbox,color=blue,vlen=0.25]
\qwire[cod=s-1-idempbox,color=orange,vlen=0.25]
\qwire[cod=s-3-idempbox,color=orange,vlen=0.25]
\qwire[cod=s-2-idempbox,color=blue,vlen=0.25]
\qwire[cod=S-outwire,vlen=0.4,color=blue]
\qwire[dom={C-u}, cod={S-3-idempbox}, color=orange]
\qwire[cod={C-cou}, dom={N-3-idempbox}, color=orange]

	}
=
	\q[scale=0.3]{
		\qwire[color=orange] \qspace \qwire[color=blue,id=outwire]\n
		\qwire[color=orange] \qcapraw[frame color=blue,hlen=3]\n
		\qwire[color=orange] \qwire[color=blue]\qcounit[color=orange,
		frame color=orange, id=cou] \qwire[color=blue]
		\n
		\qwire[color=orange]
		\qwire[color=blue]
		\qunit[color=orange,frame color=orange, id=u]
		\qwire[color=blue]
		\n
		\qwire[color=orange] \qcupraw[frame color=blue, hlen=3]\n
		\qwire[color=orange] \qspace \qwire[color=blue,vlen=1.4]
\qwire[cod=S-outwire,vlen=0.4,color=blue]
\qwire[dom={C-u}, cod={N-u}, color=orange]
\qwire[cod={C-cou}, dom={S-cou}, color=orange]

}-\frac{1}{\partial_2(\alpha_1)}
	\q[scale=0.2]{
		\qwire[color=orange] \qspace \qwire[color=blue,id=outwire]\n
		\qwire[color=orange] \qcapraw[frame color=blue,hlen=3]\n
		\qwire[color=orange] \qwire[color=blue]\qcounit[color=orange,
		frame color=orange, id=cou] \qwire[color=blue]
		\n
		\qwire[color=orange] \qunit[color=blue, frame color=blue, id=u2]
		\qwire[color=orange] \qwire[color=blue] \n
		\qcupraw[frame color=orange,hlen=3] \qwire[color=blue]
		\n
		\qspace \qwire[color=orange, id=midwire] \qspace \qwire[color=blue]\n
		\qcapraw[frame color=orange, hlen=3] \qwire[color=blue]\n
		\qwire[color=orange] \qcounit[color=blue, frame color=blue, id=cou2]
		\qwire[color=orange] \qwire[color=blue]
		\n
		\qwire[color=orange]
		\qwire[color=blue]
		\qunit[color=orange,frame color=orange, id=u]
		\qwire[color=blue]
		\n
		\qwire[color=orange] \qcupraw[frame color=blue, hlen=3]\n
		\qwire[color=orange] \qspace \qwire[color=blue,vlen=1.4]
\qwire[cod=S-outwire,vlen=0.4,color=blue]
\qwire[dom={C-u}, cod={N-u}, color=orange]
\qwire[dom={C-u2}, cod={N-u2}, color=blue]
\qwire[cod={C-cou}, dom={S-cou}, color=orange]
\qwire[cod={C-cou2}, dom={S-cou2}, color=blue]
\qwire[dom=N-midwire, color=orange, vlen=0.4]
\qwire[cod=S-midwire, color=orange, vlen=0.4]
}=\left(\partial_1(\alpha_2)-\frac{1}{\partial_2(\alpha_1)}\right)\id_{B_2B_1}
\]
and $\partial_1(\alpha_2)=2$ while
$\partial_2(\alpha_1)=1$, and so the entire three-term complex is
eliminated by $B_1$\footnote{The result continues to hold in the geometric
realization for type $\mathbb{B}_{n+1}$, as in that case,
$\partial_1(\alpha_2)=\frac{1}{\sqrt{2}}$,
which does not equal $\frac{1}{\partial_2(\alpha_1)}=\sqrt{2}$. One must
however note that the field we define our realization over can then no longer
be $\mathbb{Q}$.}.
\item  Tensor the complex on the right with $B_2$. The summand of
$P_nB_{n+1}\dots B_3B_2$  in the middle term, as well as the last two terms
in the complex, will all eliminate together in the same fashion
as~(\ref{eqn:3termelimination1}).
\\ \\
After applying equation~\ref{eqn:b2123downstepiso},
the first term and the $q^{1/2}$-shifted copy of 
$P_nB_{n+1}\dots B_3B_{212}$ will eliminate, while the remaining two terms will
eliminate since composing with the inclusion/projection morphisms of the terms
yields
\[
	\q[scale=0.5]{
	
		\qbox[hlen=3,N={},S={},n=3,s=3,id=idempbox3]
	\n
	\qspace \qspace \qcounit[color=orange, frame color=orange, id=cou]
	\n
		\qbox[hlen=3,N={},S={},n=3,s=3,id=idempbox2]
		\qwire[color=orange, id=middlewire]
	\n
	\qspace \qspace \qwire[color=orange, id=bottomwire]
	\n
	\qbox[hlen=3,N={},S={},n=3,s=3,id=idempbox]
\qwire[dom=n-3-idempbox, cod=S-bottomwire, color=orange]
\qwire[dom=C-bottomwire, cod=S-middlewire, color=orange]
\qwire[dom=N-middlewire,cod=s-3-idempbox3,color=orange]
\qwire[dom=N-bottomwire,cod=s-3-idempbox2, color=orange]
\qwire[dom=n-1-idempbox, cod=s-1-idempbox2, color=orange]
\qwire[dom=n-2-idempbox, cod=s-2-idempbox2, color=blue]
\qwire[dom=n-1-idempbox2, cod=s-1-idempbox3, color=orange]
\qwire[dom=n-2-idempbox2, cod=s-2-idempbox3, color=blue]
\qwire[dom=C-cou, cod=n-3-idempbox2,color=orange]
\qwire[cod=s-1-idempbox, color=orange,vlen=0.25]
\qwire[cod=s-2-idempbox, color=blue,vlen=0.25]
\qwire[cod=s-3-idempbox, color=orange,vlen=0.25]
\qwire[dom=n-1-idempbox3, color=orange,vlen=0.25]
\qwire[dom=n-2-idempbox3, color=blue,vlen=0.25]
\qwire[dom=n-3-idempbox3, color=orange,vlen=0.25]
}=
\q[scale=0.5]{
	
		\qbox[hlen=3,N={},S={},n=3,s=3,id=idempbox3]
	\n
	\qspace \qspace \qcounit[color=orange, frame color=orange, id=cou]
	\n
	\qwire[color=orange, id=1-idempbox2]
	\qwire[color=blue, id=2-idempbox2]
	\qwire[color=orange, id=3-idempbox2]
		\qwire[color=orange, id=middlewire]
	\n
	\qspace \qspace \qwire[color=orange, id=bottomwire]
	\n
	\qbox[hlen=3,N={},S={},n=3,s=3,id=idempbox]
\qwire[dom=n-3-idempbox, cod=S-bottomwire, color=orange]
\qwire[dom=C-bottomwire, cod=S-middlewire, color=orange]
\qwire[dom=N-middlewire,cod=s-3-idempbox3,color=orange]
\qwire[dom=N-bottomwire,cod=S-3-idempbox2, color=orange]
\qwire[dom=n-1-idempbox, cod=S-1-idempbox2, color=orange]
\qwire[dom=n-2-idempbox, cod=S-2-idempbox2, color=blue]
\qwire[dom=n-1-idempbox2, cod=s-1-idempbox3, color=orange]
\qwire[dom=n-2-idempbox2, cod=s-2-idempbox3, color=blue]
\qwire[dom=C-cou, cod=N-3-idempbox2,color=orange]
\qwire[cod=s-1-idempbox, color=orange,vlen=0.25]
\qwire[cod=s-2-idempbox, color=blue,vlen=0.25]
\qwire[cod=s-3-idempbox, color=orange,vlen=0.25]
\qwire[dom=n-1-idempbox3, color=orange,vlen=0.25]
\qwire[dom=n-2-idempbox3, color=blue,vlen=0.25]
\qwire[dom=n-3-idempbox3, color=orange,vlen=0.25]
}-
\frac{1}{\partial_2(\alpha_1)}
\q[scale=0.4]{
	
		\qbox[hlen=3,N={},S={},n=3,s=3,id=idempbox3]
	\n
	\qwire[color=orange] 
	\qunit[frame color=blue, color=blue, id=N-2-idempbox2]
	\qcounit[color=orange, frame color=orange, id=cou]
	\n
	\qcupraw[frame color=orange,hlen=3]
		\qwire[color=orange, id=middlewire]
	\n
	\qspace \qwire[color=orange, id=projectorwire]
	\n
	\qcapraw[frame color=orange,hlen=3]
	\n
	\qwire[color=orange]  \qcounit[frame color=blue, color=blue,
	id=S-2-idempbox2] \qwire[color=orange, id=bottomwire]
	\n
	\qbox[hlen=3,N={},S={},n=3,s=3,id=idempbox]
\qwire[dom=N-projectorwire, color=orange, vlen=0.4]
\qwire[cod=S-projectorwire, color=orange, vlen=0.4]
\qwire[dom=n-3-idempbox, cod=S-bottomwire, color=orange]
\qwire[dom=C-bottomwire, cod=S-middlewire, color=orange]
\qwire[dom=N-middlewire,cod=s-3-idempbox3,color=orange]
\qwire[dom=n-2-idempbox, cod=C-S-2-idempbox2, color=blue]
\qwire[dom=C-N-2-idempbox2, cod=s-2-idempbox3, color=blue]
\qwire[dom=C-cou, cod=S-cou,color=orange]
\qwire[cod=s-1-idempbox, color=orange,vlen=0.25]
\qwire[cod=s-2-idempbox, color=blue,vlen=0.25]
\qwire[cod=s-3-idempbox, color=orange,vlen=0.25]
\qwire[dom=n-1-idempbox, color=orange,vlen=0.25]
\qwire[dom=n-1-idempbox3, color=orange,vlen=0.25]
\qwire[dom=n-2-idempbox3, color=blue,vlen=0.25]
\qwire[dom=n-3-idempbox3, color=orange,vlen=0.25]
\qwire[cod=s-1-idempbox3, color=orange,vlen=0.25]
}
=
\q[scale=0.5]{\qbox[id=idempbox,hlen=3,N={},S={},n=3,s=3]
\qwire[cod=s-1-idempbox, color=orange,vlen=0.25]
\qwire[cod=s-2-idempbox, color=blue,vlen=0.25]
\qwire[cod=s-3-idempbox, color=orange,vlen=0.25]
\qwire[dom=n-1-idempbox, color=orange,vlen=0.25]
\qwire[dom=n-2-idempbox, color=blue,vlen=0.25]
\qwire[dom=n-3-idempbox, color=orange,vlen=0.25]
}
\]
as the subtracted term will become zero after applying the dot-absorption,
Frobenius, and pitchfork-killing relations.
	\end{enumerate}
The reversed versions of these eliminations all hold if one applies the same
arguments while flipping the diagrams vertically.
\end{proof}
\subsubsection{A simplified form for $P_n\JM_{n+1}$}
Given the relations we have now proven, we produce a simplified form for the
complex $P_n\JM_{n+1}$, which will ultimately give us the desired eigenmaps.
\begin{lemma}
	\label{lem:pjmsimplified}
	If $P$ is a complex in $\Ch^-(\SBim_{\mathbb{B}_{n+1}})$, where
	$n \geq 2$, such that for
	any Bott-Samelson $B_i$ where $1\leq i \leq n$, $PB_i\simeq B_iP\simeq
	0$, then the complex $P\JM_{n+1}$ is homotopy equivalent to a complex of
	the form below,
\[\begin{tikzcd}[ampersand replacement=\&,cramped, column sep=tiny]
	{q^{1/2}PB_{n+1}} \& \dots \& {q^{1/2}t^{n/2}PB_{n+1}B_n\dots B_2B_1} \\
	{q^{1/2}t^{(n+1)/2}PB_{n+1}B_n\dots B_{212}} \& \dots \& {q^{1/2}t^{n}PB_{n+1}B_n\dots B_{212}\dots B_{n}B_{n+1}} \\
	\\
	{t^{(2n+1)/2}PB_{n+1}B_n\dots B_{212}\dots B_{n}B_{n+1}} \& \dots \& {t^{3n/2}PB_{n+1}B_n\dots B_{212}} \\
	{t^{(3n+1)/2}PB_{n+1}B_n\dots B_1} \& \dots \& {t^{2n+1}P}
	\arrow[from=1-1, to=1-2]
	\arrow[from=1-2, to=1-3]
	\arrow[from=1-3, to=2-1]
	\arrow[from=2-1, to=2-2]
	\arrow[from=2-2, to=2-3]
	\arrow["\xi_n"{description}, from=2-3, to=4-1]
	\arrow[from=4-1, to=4-2]
	\arrow[from=4-2, to=4-3]
	\arrow[from=4-3, to=5-1]
	\arrow[from=5-1, to=5-2]
	\arrow[from=5-2, to=5-3]
\end{tikzcd}\]
	where all pictured arrows unless otherwise indicated are
	$\pm 1$ times the obvious
	unit or counit maps, and $\xi_n$ is
	\[
\id_{PB_{n+1}\dots B_{212}\dots B_n}\otimes \left(\mathbf{\big\vert}
	 \alpha_{n+1}-\alpha_{n+1}\mathbf{\big\vert}\right)
	\]
	such that the horizontal lines are
	identity strands for $B_{n+1}$.
\end{lemma}
\begin{proof}
	Recalling the form of $\JM_{n+1}$ from (\ref{eqn:typebjm}), we start by
	calculating the product
	\[PF(s_{n+1})F(s_n)\dots F(s_3)F(s_2)F(s_1).\]
	As $F(s_{i})=B_i\to t^{1/2}\mathbb{1}$ (where the map is the counit map),
	and $PF(s_{n+1})=PB_{n+1}\to t^{1/2}P$, we can apply
	(\ref{eqn:zerorelation1}) to see that
\begin{align*}
\begin{split}
PF(s_{n+1})\dots F(s_{i+1})B_i &\simeq PB_{n+1}\dots B_{i+1}B_i
\to \dots \to t^{(n+1-i)/2} PB_i \\
					 &\simeq PB_{n+1}\dots B_{i+1}B_i
\end{split}
\end{align*}
and therefore
\begin{align*}
\begin{split}
	PF(s_{n+1})\dots F(s_{i+1})(B_i \to t^{1/2}\mathbb{1}) \simeq
	&PB_{n+1}\dots B_{i+1} B_i \\
	&{\big\downarrow} \\
	t^{1/2}&PB_{n+1}\dots B_{i+1} \\
	&{\big\downarrow} \\
	       &\vdots\\
	&{\big\downarrow} \\
	t^{(n+1-i)/2} &PB_{n+1} \\
	       &{\big\downarrow}\\
	t^{(n+2-i)/2}&P
\end{split}
\end{align*}
giving a homotopy equivalence
\begin{equation}
\label{eqn:firstquarter}
PF(s_{n+1})\dots F(s_1)\simeq
PB_{n+1}\dots B_1 \to \dots \to t^{n/2}PB_{n+1} \to t^{(n+1)/2} P.
\end{equation}
The remaining part of $\JM_{n+1}$ is now \[F(s_2)\dots F(s_{n+1})
F(s_{n+1})\dots F(s_2)F(s_1)F(s_2)\dots F(s_{n+1}).\]
When we multiply (\ref{eqn:firstquarter}) by $F(s_2)$, we first see that
multiplying by $B_2$ and applying (\ref{eqn:zerorelation1}) and
(\ref{eqn:b2inclusion}) gives the complex
\begin{align*}
	\begin{split}
		PB_{n+1}\dots &B_{212}  \\
		&\oplus \\
		PB_{n+1}\dots &B_{2} \to t^{1/2}(q^{1/2}\mathbb{1} \oplus
		q^{-1/2}\mathbb{1})
		PB_{n+1}\dots B_{2} \to tPB_{n+1}\dots B_{2}
\end{split}
\end{align*}
where the bottom row is eliminated via (\ref{eqn:3termelimination2}), and so
the complex after multiplying by $F(s_{2})$ becomes
\begin{align*}
	PB_{n+1}\dots &B_{212}\\
		      &{\big\downarrow}\\
	t^{1/2}PB_{n+1}\dots & B_{2}B_1
	\to \dots \to t^{(n+1)/2}PB_{n+1} \to t^{(n+2)/2}P.
\end{align*}
Note that for any $n+1 \geq i > 2$,
the part of the above complex on the second row will
be eliminated by tensoring on the right with $B_i$ due to a combination of
(\ref{eqn:zerorelation1}) and (\ref{eqn:3termelimination1}), so tensoring it
with $F(s_3)\dots F(s_{n+1})$ will give us
\begin{equation}
\label{eqn:secondquarter}
\begin{tikzcd}[ampersand replacement=\&,cramped]
	{PB_{n+1}B_n\dots B_{212}\dots B_{n}B_{n+1}} \& \dots \& {t^{(n-1)/2}PB_{n+1}B_n\dots B_{212}} \\
	{t^{n/2}PB_{n+1}B_n\dots B_1} \& \dots \& {t^{(2n+1)/2}P}
	\arrow[from=1-1, to=1-2]
	\arrow[from=1-2, to=1-3]
	\arrow[from=1-3, to=2-1]
	\arrow[from=2-1, to=2-2]
	\arrow[from=2-2, to=2-3]
\end{tikzcd}.
\end{equation}
When we multiply by $F(s_{n+1})$ again, we note that multiplication by
$B_{n+1}$ will contract the entire second row of the above complex via
(\ref{eqn:zerorelation1}) and (\ref{eqn:3termelimination1}), and eliminate all
terms on the first row except
\begin{align*}
	\begin{split}
		q^{1/2}PB_{n+1}\dots &B_{212} \dots B_{n+1}\\
		&\oplus\\
		q^{-1/2}PB_{n+1}\dots &B_{212}\dots B_{n+1} \to t^{1/2}
		PB_{n+1}\dots B_{212} \dots B_{n+1}
	\end{split}
\end{align*}
where the second row will be eliminated by being an instance of the lower arrow
in (\ref{eqn:3termelimination1}). After adding the compensatory term from this
elimination to the arrow given by $B_{n+1}\to t^{1/2}\mathbb{1}$, the complex
is now
\[
\begin{tikzcd}[ampersand replacement=\&]
	{q^{1/2}PB_{n+1}B_n\dots B_{212}\dots B_{n}B_{n+1}} \&\& \\
	{t^{1/2}PB_{n+1}B_n\dots B_{212}\dots B_{n}B_{n+1}} \& \dots \& {t^{n/2}PB_{n+1}B_n\dots B_{212}} \\
	{t^{(n+1)/2}PB_{n+1}B_n\dots B_1} \& \dots \& {t^{(2n+2)/2}P}
	\arrow["\xi_n"
	,shift right=20, from=1-1, to=2-1]
	\arrow[from=2-1, to=2-2]
	\arrow[from=2-2, to=2-3]
	\arrow[from=2-3, to=3-1]
	\arrow[from=3-1, to=3-2]
	\arrow[from=3-2, to=3-3]
\end{tikzcd}\]
where $\xi_n$ is as described in the statement of the lemma.
When multiplying by
the remaining part
$F(s_{n})\dots F(s_1)\dots F(s_{n})F(s_{n+1})$ of
$\JM_{n+1}$, we note that for $i < n+1$, $B_i$ will eliminate the entire part of the
complex without a $q$-shift via (\ref{eqn:zerorelation1})
and the subcomplex eliminations in
section~\ref{ssec:subcomplexeliminations}, and eliminate all but the lowest
$t$-shifted term in the $q^{1/2}$-shifted part of the complex.
\\ \\
Then, composing
the isomorphism (\ref{eqn:downshiftiso})
(or (\ref{eqn:b21downstepiso}) when doing the same to a term with
$B_{212}$ as the rightmost factor) 
with the counit gives that the arrow from the new lowest-$t$-degree term is
now a unit.
Thus, right before multiplying by $F(s_{n+1})$, the complex
is homotopy equivalent to
\[\begin{tikzcd}[ampersand replacement=\&,cramped,column sep=tiny]
	{q^{1/2}PB_{n+1}B_{n}} \& \dots \& {q^{1/2}t^{(n-1)/2}PB_{n+1}B_n\dots B_2B_1} \\
	{q^{1/2}t^{n/2}PB_{n+1}B_n\dots B_{212}} \& \dots \& {q^{1/2}t^{(2n-1)/2}PB_{n+1}B_n\dots B_{212}\dots B_{n}B_{n+1}} \\
	\\
	{t^{2n/2}PB_{n+1}B_n\dots B_{212}\dots B_{n}B_{n+1}} \& \dots \& {t^{(3n-1)/2}PB_{n+1}B_n\dots B_{212}} \\
	{t^{3n/2}PB_{n+1}B_n\dots B_1} \& \dots \& {t^{(4n+1)/2}P}
	\arrow[from=1-1, to=1-2]
	\arrow[from=1-2, to=1-3]
	\arrow[from=1-3, to=2-1]
	\arrow[from=2-1, to=2-2]
	\arrow[from=2-2, to=2-3]
	\arrow["\xi_n"{description}, from=2-3, to=4-1]
	\arrow[from=4-1, to=4-2]
	\arrow[from=4-2, to=4-3]
	\arrow[from=4-3, to=5-1]
	\arrow[from=5-1, to=5-2]
	\arrow[from=5-2, to=5-3]
\end{tikzcd}\]
where all the unlabeled arrows in the $q$-shifted part are units and all
unlabeled arrows in the part without a $q$-shift are counits. To complete
the proof, it remains to show that $B_{n+1}$ (which obviously eliminates the
rest of the complex besides the lowest-$t$-degree term and
the middle terms) will also eliminate the
middle two terms at degree $t^{(2n-1)}$ and $t^{2n/2}$.
\\ \\
After eliminating the
$q$-shifted and $q^{-1/2}$-shifted terms via subcomplex eliminations, the remaining
two terms eliminate by virtue of the fact that
\begin{equation}
	\label{eqn:finaleliminationmap}
\q[scale=0.6]{\qspider[frame color=none,
n=3,s=3,N={},S={},id=mul]
\qwire[cod=N-2-mul,dom=C-mul]
\qwire[dom=S-1-mul,cod=C-mul]
\qwire[dom=S-3-mul,cod=C-mul] \n
\qwire \qbox[frame color=none, n=0,s=0,hlen=1, name=$\alpha_{n+1}$] \qwire \n
\qspider[frame color=none,
n=3,s=3,N={},S={},id=comul]
\qwire[cod=S-2-comul,dom=C-comul]
\qwire[dom=N-1-comul,cod=C-comul]
\qwire[dom=N-3-comul,cod=C-comul]
}-\q[scale=0.6]{\qspider[frame color=none,
n=3,s=3,N={},S={},id=mul]
\qwire[cod=N-2-mul,dom=C-mul]
\qwire[dom=S-1-mul,cod=C-mul]
\qwire[dom=S-3-mul,cod=C-mul] \n
\qwire \qspace \qwire \qbox[n=0,s=0,
hlen=1,frame color=none, name=$\alpha_{n+1}$]\n
\qspider[frame color=none,
n=3,s=3,N={},S={},id=comul]
\qwire[cod=S-2-comul,dom=C-comul]
\qwire[dom=N-1-comul,cod=C-comul]
\qwire[dom=N-3-comul,cod=C-comul]
}=2\q[scale=0.4]{
\qspider[frame color=none,
n=3,s=3,N={},S={},id=mul]
\qwire[cod=N-2-mul,dom=C-mul]
\qwire[dom=S-1-mul,cod=C-mul]
\qwire[dom=S-3-mul,cod=C-mul] \n
\qwire \qspace \qunit[color=black] \n
\qwire \qspace \qcounit[color=black]\n
\qspider[frame color=none,
n=3,s=3,N={},S={},id=comul]
\qwire[cod=S-2-comul,dom=C-comul]
\qwire[dom=N-1-comul,cod=C-comul]
\qwire[dom=N-3-comul,cod=C-comul]
}-0=2\id_{B_{n+1}} \end{equation}
after composing the inclusion/projection maps with
$\xi_n$. The complex
$P\JM_{n+1}$ is thus of the form given in the statement of the lemma. Since all
eliminations were done on copies of the complex multiplied by a non-identity
Bott-Samelson, the pictured arrows all continue to be the expected unit/counit
maps.
\end{proof}
\noindent We can now prove statements about the complex of morphisms
$\Hom(P_n,P_n\JM_{n+1})$.
\subsubsection{Existence of the eigenmaps}
Once we apply the $\Hom(P_n,-)$ functor
(equivalent to the $\Hom(P_n,P_n\otimes -)$ functor and therefore
$\Hom(\mathbb{1},-)$
as all terms have a tensor factor of $P_n$ on the left)
to the entire complex given by
Lemma~\ref{lem:pjmsimplified}, it is immediately obvious, by the fact that
$\SBim_{\mathbb{B}_{n+1}}$ is a rigid monoidal category and $P_n$ kills lower Bott-Samelsons,
that all terms besides those involving $P$, $PB_{n+1}$, or 
$PB_{n+1}\dots B_{212}\dots B_{n+1}$ are taken to zero
via application
of Theorem~\ref{thm:idempotentmorphismzero}. We now show that the
latter of these types of terms is also taken to zero.
\begin{lemma}
	\label{lem:homspacezerolightleaves}
	If $P$ is a complex in $\Ch^-(\SBim_{\mathbb{B}_{n+1}})$ such that for
	any Bott-Samelson $B_i$ where $1 \leq i \leq n$, $PB_i\simeq
	B_iP=0$, then
	\[\Hom
		(\mathbb{1},PB_{n+1}B_nB_{n-1}\dots B_{3}B_{212}B_{3}\dots
	B_{n-1}B_nB_{n+1})\simeq 0.\]
\end{lemma}
\begin{proof}
It is sufficient to show that all morphisms in the complex are linear combinations
of morphisms that factor through a product of $P$ with a
non-identity Bott-Samelson in
$\SBim_{\mathbb{B}_{n}}$ and are therefore
$0$. A morphism from $\mathbb{1}$ to a component of $PB_{n+1}\dots B_{212}\dots
B_{n+1}$, factoring through a morphism from $\mathbb{1}$ to
$XB_{n+1}\dots B_{212}\dots B_{n+1}$ for some
Bott-Samelson bimodule
$X\in \SBim_{\mathbb{B}_{n+1}}$ corresponding to a word
$\underline{w}$ in $\mathbb{B}_n$, will be an $R$-linear combination of
morphisms given by the double leaves basis~\cite{elias_introduction_2020}, indexed
by sub-expressions of $\underline{w}
s_{n+1}s_n \dots s_2 s_1s_2
\dots s_{n} s_{n+1}$ that are the identity, post-composed with the
projector for $B_{212}$ at the appropriate strands.
\\ \\
Since $\underline{w}$ does not contain $s_{n+1}$, any sub-expression
of $\underline{w}s_{n+1}\dots s_{2}s_1s_2\dots
s_{n+1}$ equivalent to the identity must either include both copies of
$s_{n+1}$ or exclude both of them. Suppose both are excluded. Then, applying
the double leaves basis algorithm will leave us with a morphism that factors
through $PB_{n}\dots B_{212}\dots B_n$, which will therefore be zero.

\[\q[scale=0.8]{

	\qbox[color=none,frame color=none,n=0,s=0,name=$\cdots$,id=ellipsislast]
	\qunit[color=black]\symbolO[1]{$s_{n+1}$}
	\qwire[color=red]\symbolO[1]{$s_n$}
	\qbox[color=none,frame color=none,n=0,s=0,name=$\cdots$,id=ellipsisleft]

		\qwire[color=green,vlen=1,id=3wireleft]\symbolO[1]{$s_3$}
		\qbox[n=3,s=3,N={},S={},hlen=3,id=idempbox,name=$B_{212}$]
		\symbolN[1]{$s_2$}[id=idempbox]
		\symbolN[3]{$s_2$}[id=idempbox]
		\symbolN[2]{$s_1$}[id=idempbox]
		\qwire[color=green,vlen=1,id=3wireright]\symbolO[1]{$s_3$}


	\qbox[color=none,frame color=none,n=0,s=0,name=$\cdots$,id=ellipsisright]

	\qwire[color=red]\symbolO[1]{$s_n$}

	\qunit[color=black]\symbolO[1]{$s_{n+1}$}
\n
\qspace \qspace \qwire[color=red] \qspace \qwire[color=green]
\qwire[color=orange]\qwire[color=blue] \qwire[color=orange] \qwire[color=green]
\qspace \qwire[color=red]
\n
\qspace\qspace
\qbox[color=none,frame color=none,n=0,s=0,name=$\cdots$]

\qbox[color=none,frame color=none,n=0,s=0,name=$\cdots$]
\qbox[color=none,frame color=none,n=0,s=0,name=$\cdots$]
\qbox[color=none,frame color=none,n=0,s=0,name=$\cdots$]
\qbox[color=none,frame color=none,n=0,s=0,name=$\cdots$]
\qbox[color=none,frame color=none,n=0,s=0,name=$\cdots$]
\qbox[color=none,frame color=none,n=0,s=0,name=$\cdots$]
\qbox[color=none,frame color=none,n=0,s=0,name=$\cdots$]
\qbox[color=none,frame color=none,n=0,s=0,name=$\cdots$]

\qwire[dom=n-1-idempbox,color=orange,vlen=0.25]
\qwire[dom=n-3-idempbox,color=orange,vlen=0.25]
\qwire[dom=n-2-idempbox,color=blue,vlen=0.25]
\qwire[cod=s-1-idempbox,color=orange,vlen=0.25]
\qwire[cod=s-3-idempbox,color=orange,vlen=0.25]
\qwire[cod=s-2-idempbox,color=blue,vlen=0.25]
}\]
Now, if both $s_{n+1}$ are included in the sub-expression, they must be
joined by a line in the corresponding morphism. If the first $s_n$
from the right is
included in the sub-expression, the non-commuting of $s_n$ and $s_{n+1}$
means that it can only be cancelled out by the second $s_n$ on the right, so,
once again, either both $s_n$ are included or both are excluded. If both
are excluded, then we once again have a morphism that is zero, since it
factors through $PB_{n-1}\dots B_{212}\dots B_{n-1}$.
\[\q[scale=0.8]{

	\qbox[color=none,frame color=none,n=0,s=0,name=$\cdots$,id=ellipsislast]
	\qwire[id=nwireleft]\symbolO[1]{$s_{n+1}$}

	\qwire[color=red, id=ithwireleft]\symbolO[1]{$s_n$}
	\qbox[color=none,frame color=none,n=0,s=0,name=$\cdots$]
		\qwire[color=green,vlen=1,id=3wireleft]\symbolO[1]{$s_3$}
		\qbox[n=3,s=3,N={},S={},hlen=3,id=idempbox,name=$B_{212}$]
		\symbolN[1]{$s_2$}[id=idempbox]
		\symbolN[3]{$s_2$}[id=idempbox]
		\symbolN[2]{$s_1$}[id=idempbox]
		\qwire[color=green,vlen=1,id=3wireright]\symbolO[1]{$s_3$}

	\qbox[color=none,frame color=none,n=0,s=0,name=$\cdots$]
	\qwire[color=red, id=ithwireright]\symbolO[1]{$s_n$}

	\qwire[id=nwireright]\symbolO[1]{$s_{n+1}$}
\n
\qspace\qwire \qunit[color=red,frame color=red, id=nunitleft]
\qspace \qspace\qspace\qspace
\qspace\qspace\qspace \qunit[color=red, frame color=red, id=nunitright]
\qwire
\n
\qspace\qcupraw[n=2, hlen=11]
\n\n\n\n
\qspace \qspace \qspace
\qbox[color=none,frame color=none,n=0,s=0,name=$\cdots$]
\qbox[color=none,frame color=none,n=0,s=0,name=$\cdots$, id=bottom3left
]
\qbox[color=none,frame color=none,n=1,N={},s=0,name=$\cdots$, id=bottom2left]
\qbox[color=none,frame color=none,n=1,N={},s=0,name=$\cdots$, id=bottom1]
\qbox[color=none,frame color=none,n=1,N={},s=0,name=$\cdots$, id=bottom2right]
\qbox[color=none,frame color=none,n=1,N={},s=0,name=$\cdots$, id=bottom3right]
\qbox[color=none,frame color=none,n=0,s=0,name=$\cdots$]

\qwire[dom={n-1-bottom3left},cod={I-3wireleft},color=green]
\qwire[dom={n-1-bottom3right},cod={I-3wireright},color=green]
\qwire[dom={n-1-bottom2left},cod={s-1-idempbox},color=orange]
\qwire[dom={n-1-bottom2right},cod={s-3-idempbox},color=orange]
\qwire[dom={n-1-bottom1},cod={s-2-idempbox},color=blue]

\qwire[dom=n-1-idempbox,color=orange,vlen=0.25]
\qwire[dom=n-3-idempbox,color=orange,vlen=0.25]
\qwire[dom=n-2-idempbox,color=blue,vlen=0.25]
\qwire[cod=s-1-idempbox,color=orange,vlen=0.25]
\qwire[cod=s-3-idempbox,color=orange,vlen=0.25]
\qwire[cod=s-2-idempbox,color=blue,vlen=0.25]
\qwire[dom={C-nunitleft}, cod={N-1-nunitleft}, color=red]
\qwire[dom={C-nunitright}, cod={N-1-nunitright}, color=red]
}\]
If both $s_n$ are included, then we apply the same argument to the two copies
of
$s_{n-1}$ in the expression, and find the morphism becomes zero by the
same commutativity argument if both are excluded. We can continue this recursively
for any expression where the pairs of $s_{n+1},s_{n},\dots, s_{i+1}$
for $n+1 \geq i > 2$
are all included in the sub-expression, and split into a case where
both copies of $s_i$ are excluded (making the morphism zero), or where both
are included (which splits into the cases where we include or exclude
both copies of $s_{i-1}$).
\[
	\q[scale=0.8]{

	\qbox[color=none,frame color=none,n=0,s=0,name=$\cdots$,id=ellipsislast]
	\qwire[id=nwireleft]\symbolO[1]{$s_{n+1}$}

	\qwire[color=red]\symbolO[1]{$s_n$}
	\qbox[color=none,frame color=none,n=0,s=0,name=$\cdots$]
	\qwire[color=purple, id=ithwireleft]\symbolO[1]{$s_i$}
	\qbox[color=none,frame color=none,n=0,s=0,name=$\cdots$]
		\qwire[color=green,vlen=1,id=3wireleft]\symbolO[1]{$s_3$}
		\qbox[n=3,s=3,N={},S={},hlen=3,id=idempbox,name=$B_{212}$]
		\symbolN[1]{$s_2$}[id=idempbox]
		\symbolN[3]{$s_2$}[id=idempbox]
		\symbolN[2]{$s_1$}[id=idempbox]
		\qwire[color=green,vlen=1,id=3wireright]\symbolO[1]{$s_3$}

\qbox[color=none,frame color=none,n=0,s=0,name=$\cdots$]
\qwire[color=purple, id=ithwireright]\symbolO[1]{$s_i$}

	\qbox[color=none,frame color=none,n=0,s=0,name=$\cdots$]
	\qwire[color=red, id=ithwireright]\symbolO[1]{$s_n$}

	\qwire[id=nwireright]\symbolO[1]{$s_{n+1}$}
\n
\qspace \qwire \qwire[color=red]
\qspace
\qunit[color=purple, frame color=purple, id=ithunitleft]
\qspace \qspace \qspace \qspace \qspace \qspace \qspace

\qunit[color=purple, frame color=purple, id=ithunitright]
\qspace
	\qwire[color=red]
	\qwire
\n
\qspace\qwire \qcupraw[n=2,hlen=13, frame color=red] \qwire
\n
\qspace\qcupraw[n=2, hlen=15]
\n\n\n\n\n
\qspace\qspace\qspace\qspace\qspace
\qbox[color=none,frame color=none,n=0,s=0,name=$\cdots$]
\qbox[color=none,frame color=none,n=0,s=0,name=$\cdots$, id=bottom3left
]
\qbox[color=none,frame color=none,n=1,N={},s=0,name=$\cdots$, id=bottom2left]
\qbox[color=none,frame color=none,n=1,N={},s=0,name=$\cdots$, id=bottom1]
\qbox[color=none,frame color=none,n=1,N={},s=0,name=$\cdots$, id=bottom2right]
\qbox[color=none,frame color=none,n=1,N={},s=0,name=$\cdots$, id=bottom3right]
\qbox[color=none,frame color=none,n=0,s=0,name=$\cdots$]
\qwire[dom={n-1-bottom3left},cod={I-3wireleft},color=green]
\qwire[dom={n-1-bottom3right},cod={I-3wireright},color=green]
\qwire[dom={n-1-bottom2left},cod={s-1-idempbox},color=orange]
\qwire[dom={n-1-bottom2right},cod={s-3-idempbox},color=orange]
\qwire[dom={n-1-bottom1},cod={s-2-idempbox},color=blue]

\qwire[dom=n-1-idempbox,color=orange,vlen=0.25]
\qwire[dom=n-3-idempbox,color=orange,vlen=0.25]
\qwire[dom=n-2-idempbox,color=blue,vlen=0.25]
\qwire[cod=s-1-idempbox,color=orange,vlen=0.25]
\qwire[cod=s-3-idempbox,color=orange,vlen=0.25]
\qwire[cod=s-2-idempbox,color=blue,vlen=0.25]
\qwire[dom={C-ithunitright}, cod={N-1-ithunitright}, color=purple]
\qwire[dom={C-ithunitleft}, cod={N-1-ithunitleft}, color=purple]
}\]
When we reach the case where we must include or exclude the copies of 
$s_{2}$ from the expression, excluding both of them results in a morphism
factoring through $PB_1$, which is zero, and including both of them results
in a morphism that will be killed by $B_{212}$ via the pitchfork-killing relation.
\[\q[scale=0.8]{

	\qbox[color=none,frame color=none,n=0,s=0,name=$\cdots$,id=ellipsislast]
	\qwire[id=nwireleft]\symbolO[1]{$s_{n+1}$}
	\qbox[color=none,frame color=none,n=0,s=0,name=$\cdots$,id=ellipsisleft]

		\qwire[color=green,vlen=1,id=3wireleft]\symbolO[1]{$s_3$}
		\qbox[n=3,s=3,N={},S={},hlen=3,id=idempbox,name=$B_{212}$]
		\symbolN[1]{$s_2$}[id=idempbox]
		\symbolN[3]{$s_2$}[id=idempbox]
		\symbolN[2]{$s_1$}[id=idempbox]
		\qwire[color=green,vlen=1,id=3wireright]\symbolO[1]{$s_3$}

	\qbox[color=none,frame color=none,n=0,s=0,name=$\cdots$,id=ellipsisright]
	\qwire[id=nwireright]\symbolO[1]{$s_{n+1}$}
\n
\qspace \qwire \qspace \qwire[color=green]
\qwire[color=orange] \qunit[color=blue, frame color=blue,
id=u] 
\qwire[color=orange]
\qwire[color=green] \qspace \qwire
\n
\qspace \qwire \qspace \qwire[color=green]
\qcupraw[n=2,hlen=3,frame color=orange] \qwire[color=green] \qspace
 \qwire
\n
\qspace \qwire
 \qspace
\qcupraw[n=2,hlen=5,frame
color=green] \qspace  \qwire
\n
\qspace\qcupraw[n=2, hlen=9]

\qwire[dom=n-1-idempbox,color=orange,vlen=0.25]
\qwire[dom=n-3-idempbox,color=orange,vlen=0.25]
\qwire[dom=n-2-idempbox,color=blue,vlen=0.25]
\qwire[cod=s-1-idempbox,color=orange,vlen=0.25]
\qwire[cod=s-3-idempbox,color=orange,vlen=0.25]
\qwire[cod=s-2-idempbox,color=blue,vlen=0.25]
\qwire[dom={C-u}, cod={N-1-u}, color=blue]
}=0\]
Therefore, all morphisms from $\mathbb{1}$ to
$PB_{n+1}\dots B_{212}\dots B_{n+1}$ have as their components morphisms that
factor through the product of $P$ with a non-identity Soergel bimodule in
$\SBim_{\mathbb{B}_n}$, which will be zero.
\end{proof}
\begin{corollary}
	\label{cor:homcomplex}
	The complex $\Hom(P_n,P_n\JM_{n+1})$ is homotopy equivalent to
\begin{align*}
	\Hom &(P_n,q^{1/2}P_nB_{n+1}) \\
	     &{\big\downarrow}\\
	 \biggl( \Hom(P_n,t^{(4n+1)/2}P_nB_{n+1})
	     & \xrightarrow{\Hom(P_n,\id_{P_n} \otimes
	     \q[scale=0.4]{\qcounit[color=black]})}
     \Hom(P_n,t^{2n+1}P_n) \biggr)
\end{align*}
\end{corollary}
\begin{proof}
	To the form of the complex given by Lemma~\ref{lem:pjmsimplified},
	apply the fact that
	$\SBim_{\mathbb{B}_{n+1}}$ is a rigid
	monoidal category and then
	apply Lemma~\ref{lem:homspacezerolightleaves}. As the twists in the complex
	were all from strictly lower-$t$-degree terms to higher $t$-degree terms,
	the arrow $\Hom(P,t^{(4n+1)/2}PB_{n}+1)\to \Hom(P,t^{2n+1}P)$ is unchanged
	by homological perturbation.
\end{proof}
\noindent To prove the existence of the eigenmap $\beta_n$, we show that
the space of homotopy classes of morphisms
$\Hom(P_n,P_n\JM_{n+1})$ surjects onto the homotopy classes of morphisms
$\Hom(P_n,q^{1/2}P_nB_{n+1})$.
\begin{lemma}
	\label{lem:homspacesurjection}
	If the complex $\End(P_n)$ has cohomology concentrated in degrees of the
	same parity, then
	there is a surjection of cohomologies
	\[H^\bullet(\Hom(P_n,P_n\JM_{n+1})) \twoheadrightarrow
		H^\bullet(\Hom(P_n,q^{1/2}P_nB_{n+1})).
	\]
\end{lemma}
\begin{proof}
	Denote the complex $\Hom(P_n,t^{(4n+1)/2}P_nB_{n+1}\to t^{2n+1}P_n)$ as $C$.
	Corollary~\ref{cor:homcomplex} gives a short exact sequence of complexes
	\[0\to C \to \Hom(P_n,P_n\JM_{n+1}) \to
	\Hom(P_n,q^{1/2}P_nB_{n+1})\]
	which induces a long exact sequence of cohomologies
	\[\begin{tikzcd}[ampersand replacement=\&,cramped]
	\dots \&\& \\
	{H^i(C)} \& {H^i(\Hom(P_n,P_n\JM_{n+1}))} \& {H^i(\Hom(P_n,q^{1/2}P_{n}B_{n+1}))} \\
	{H^{i+1}(C)} \& {H^{i+1}(\Hom(P_n,P_n\JM_{n+1}))} \& {H^{i+1}(\Hom(P_n,q^{1/2}P_{n}B_{n+1}))} \\
	\&\& \dots
	\arrow[from=1-1, to=2-1]
	\arrow[from=2-1, to=2-2]
	\arrow[from=2-2, to=2-3]
	\arrow[from=2-3, to=3-1]
	\arrow[from=3-1, to=3-2]
	\arrow[from=3-2, to=3-3]
	\arrow[from=3-3, to=4-3]
\end{tikzcd}\]
where, if both $C$ and $\Hom(P_n,q^{1/2}P_nB_{n+1})$ have cohomology
concentrated in degrees of the same parity, the desired surjection is achieved.
\\ \\
We assert that $\Hom(P_n,q^{1/2}P_nB_{n+1})$ has the same parity as
$\End(P_n)$. This follows from the fact that, by the double leaves basis
and the fact that $P_n$ does not contain any copies of $B_{n+1}$ (and so
when computing a basis element, the copy of $s_{n+1}$ in each word must
be excluded), any
morphism $P_n\to q^{1/2}P_nB_{n+1}$ is a morphism $P_n\to P_{n}$ tensored with
the unit for $B_{n+1}$. It now remains to prove that $C$ is also
parity in cohomology.
\\ \\
Observe that there is a long exact sequence of cohomologies
\[\begin{tikzcd}[ampersand replacement=\&,cramped]
	\dots \&\& \\
	{H^{i-1}(\Hom(P_n,t^{2n+1}P_n))} \& {H^{i-1}(C)} \& {H^{i-1}(\Hom(P_n,t^{(4n+1)/2}P_nB_{n+1}))} \\
	{H^i(\Hom(P_n,t^{2n+1}P_n))} \& {H^i(C)} \& {H^{i}(\Hom(P_n,t^{(4n+1)/2}P_nB_{n+1}))} \\
	{H^{i+1}(\Hom(P_n,t^{2n+1}P_n))} \& {H^{i+1}(C)} \& {H^{i+1}(\Hom(P_n,t^{(4n+1)/2}P_nB_{n+1}))} \\
	\&\& \dots
	\arrow[from=1-1, to=2-1]
	\arrow[from=2-1, to=2-2]
	\arrow[from=2-2, to=2-3]
	\arrow["\delta"{description}, from=2-3, to=3-1]
	\arrow[from=3-1, to=3-2]
	\arrow[from=3-2, to=3-3]
	\arrow["\delta"{description}, from=3-3, to=4-1]
	\arrow[from=4-1, to=4-2]
	\arrow[from=4-2, to=4-3]
	\arrow[from=4-3, to=5-3]
\end{tikzcd}\]
where $\delta$ is the map induced on cohomology by post-composition with
$\id_{P_n} \otimes \q[scale=0.4]{\qcounit[color=black]}$. Additionally
observe that, assuming $i$ is of the parity where $\End(P_n)$'s cohomology
is concentrated, and since $\Hom(P_n,t^{2n+1}P_n)$ is an even shift of
$\End(P_n)$, the long exact sequence actually looks like
\[\begin{tikzcd}[ampersand replacement=\&,cramped]
	\dots \&\& \\
	0 \& {H^{i-1}(C)} \& {H^{i-1}(\Hom(P_n,t^{(4n+1)/2}P_nB_{n+1}))} \\
	{H^i(\Hom(P_n,t^{2n+1}P_n))} \& {H^i(C)} \& 0 \\
	0 \& {H^{i+1}(C)} \& {H^{i+1}(\Hom(P_n,t^{(4n+1)/2}P_nB_{n+1}))} \\
	\&\& \dots
	\arrow[from=1-1, to=2-1]
	\arrow[from=2-1, to=2-2]
	\arrow[from=2-2, to=2-3]
	\arrow["\delta"{description}, from=2-3, to=3-1]
	\arrow[from=3-1, to=3-2]
	\arrow[from=3-2, to=3-3]
	\arrow[from=3-3, to=4-1]
	\arrow[from=4-1, to=4-2]
	\arrow[from=4-2, to=4-3]
	\arrow[from=4-3, to=5-3]
\end{tikzcd}\]
additionally, since a morphism from $P_n$ to $P_nB_{n+1}$ is simply an
endomorphism of $P_n$ tensored with the unit for $B_{n+1}$, post-composing with
the counit will mean that $\delta$ is multiplication by $\alpha_{n+1}$, which
is injective, meaning $C$ has the same parity in cohomology as
$\End(P_n)$ and therefore
$\Hom(P_n,q^{1/2}P_nB_{n+1})$.
\end{proof}
\noindent From this, it follows that there exists a cochain map
$\beta_n: qP_n\to P_n\JM_{n+1}$ of the form
\[\begin{tikzcd}[ampersand replacement=\&,cramped,column sep=tiny]
	{qP_n} \&\& \\
	{q^{1/2}P_nB_{n+1}} \& \dots \& {q^{1/2}t^{n/2}P_nB_{n+1}B_n\dots B_2B_1} \\
	{q^{1/2}t^{(n+1)/2}P_nB_{n+1}B_n\dots B_{212}} \& \dots \& {q^{1/2}t^{n}P_nB_{n+1}B_n\dots B_{212}\dots B_{n}B_{n+1}} \\
	\\
	{t^{(2n+1)/2}P_nB_{n+1}B_n\dots B_{212}\dots B_{n}B_{n+1}} \& \dots \& {t^{3n/2}P_nB_{n+1}B_n\dots B_{212}} \\
	{t^{(3n+1)/2}P_nB_{n+1}B_n\dots B_1} \& \dots \& {t^{2n+1}P_n}
	\arrow[color=red, from=1-1, to=2-1]
	\arrow[from=2-1, to=2-2]
	\arrow[from=2-2, to=2-3]
	\arrow[from=2-3, to=3-1]
	\arrow[from=3-1, to=3-2]
	\arrow[from=3-2, to=3-3]
	\arrow["{\xi_n}"{description}, from=3-3, to=5-1]
	\arrow[from=5-1, to=5-2]
	\arrow[from=5-2, to=5-3]
	\arrow[from=5-3, to=6-1]
	\arrow[from=6-1, to=6-2]
	\arrow[from=6-2, to=6-3]
\end{tikzcd}\]
where the red arrow is $\id_{P_n}\otimes \q[scale=0.4]{\qunit[color=black]}$.
We are now able to prove that $\Cone(\beta_n)$ kills non-identity Bott-Samelsons on
both the left and right.
\begin{lemma}
	\label{lem:conebeta}
	Given $\beta_n:qP_n\to P_n\JM_{n+1}$ defined as above,
	$B_{i}\Cone(\beta_n)\simeq \Cone(\beta_n)B_i\simeq 0$ for all
	non-identity Bott-Samelsons $B_i\in\SBim_{\mathbb{B}_{n+1}}$.
\end{lemma}
\begin{proof}
Clearly, $\Cone(\beta_n)B_i\simeq 0 $ via the multiplication rules and
subcomplex eliminations given in Section~\ref{ssec:subcomplexeliminations}.
We show that $B_i\Cone(\beta_n)\simeq 0$ using the
existence
and properties of the categorical idempotent in Theorem~\ref{thm:sbimidempotent}.
\\ \\
Take the map \[\Cone(\beta_n) \xrightarrow{\id_{\Cone(\beta_n)}\otimes
\eta_{\mathbb{B}_{n+1}}} \Cone(\beta_n)\mathbf{P}_{\mathbb{B}_n}.\]
The cone of this map
is given by $\Cone(\beta_n)\otimes \Cone(\eta_{\mathbb{B}_{n+1}})$, but since
$\Cone(\beta_n)$ kills non-identity Soergel bimodules on its right and
$\Cone(\eta_{\mathbb{B}_{n+1}})$ consists only of non-identity Soergel bimodules, the
cone is contractible and the map is a homotopy equivalence. Since
$\mathbf{P}_{\mathbb{B}_{n+1}}$ is in the Drinfeld center,
\[\Cone(\beta_n)\simeq \Cone(\beta_n)\mathbf{P}_{\mathbb{B}_{n+1}}\simeq
\mathbf{P}_{\mathbb{B}_{n+1}}\Cone(\beta_n)\]
which kills non-identity Bott-Samelsons on its left as well.
\end{proof}
\noindent The map $\alpha_n:t^{2n+1}P_n \to P_n\JM_{n+1}$ is given by the inclusion 
of the last copy of $P_n$ in the complex
from Lemma~\ref{lem:pjmsimplified}, and since its mapping cone can be
simplified so that the two copies of $P_n$ are eliminated, the remaining
terms are all tensored with some non-identity Soergel bimodule.
\subsubsection{Constructing $P_{n+1}'$}
\label{ssec:constructpnplusone}
We use the shorthand $\frac{q}{1-qt^{-(2n+1)}}P$ to denote the following
semi-infinite direct sum of complexes:
\[\bigoplus_{i=0}^n q(qt^{-(2n+1)})^i P_n
\]
and $\frac{1}{1-qt^{-(2n+1)}}P_n\JM_{n+1}$ to denote the semi-infinite direct
sum
\[\bigoplus_{i=0}^n (qt^{-(2n+1)})^i P_n\JM_{n+1}.
\]
We note that $\alpha_n$ and $\beta_n$ induce cochain maps
\[\frac{q}{1-qt^{-(2n+1)}}P\xrightarrow{\alpha_n',\beta_n'}
\frac{1}{1-qt^{-(2n+1)}}P_n\JM_{n+1}\]
where $\alpha_n'$ is componentwise
\[\alpha_n:q(qt^{-(2n+1)})^i P_n\to
(qt^{-(2n+1)})^{i+1}P_n\JM_{n+1}\]
and
$\beta_n'$ is componentwise
\[\beta_n:q(qt^{-(2n+1)})^i\to (qt^{-(2n+1)})^iP_n\JM_{n+1}\].
Thus, $\alpha_n'+\beta_n'$ forms a ``ladder'' which appears as follows
(taking all arrows to be degree zero):
\[\begin{tikzcd}[ampersand replacement=\&,cramped]
	{qP_{n}} \& {P_{n}\JM_{n+1}} \\
	{q^2t^{-(2n+1)}P_{n}} \& {qt^{-(2n+1)}P_{n}\JM_{n+1}} \\
	{q^3t^{-2(2n+1)}P_{n}} \& {q^2t^{-2(2n+1)}P_{n}\JM_{n+1}} \\
	\dots \& \dots
	\arrow["{\beta_n}", from=1-1, to=1-2]
	\arrow["{\alpha_n}", from=1-1, to=2-2]
	\arrow["{\beta_n}", from=2-1, to=2-2]
	\arrow["{\alpha_n}", from=2-1, to=3-2]
	\arrow["{\beta_n}", from=3-1, to=3-2]
	\arrow["{\alpha_n}", from=3-1, to=4-2]
\end{tikzcd}\]
If we write $\Cone(\alpha_n'+\beta_n')$ as an extension of complexes (therefore
adding a shift of $(tq)^{-1/2}$ to $\frac{q}{1-qt^{-(2n+1)}}P_n$), we can
use it to
construct the semi-infinite
extension of complexes
\[
	P_{n+1}':=t^{-(2n+1)}\left(
\begin{tikzcd}[ampersand replacement=\&,cramped]
	{(tq)^{-1/2}qP_{n}} \& {P_{n}\JM_{n+1}} \\
	{(tq)^{-1/2}q^2t^{-(2n+1)}P_{n}} \& {qt^{-(2n+1)}P_{n}\JM_{n+1}} \\
	{(tq)^{-1/2}q^3t^{-2(2n+1)}P_{n}} \& {q^2t^{-2(2n+1)}P_{n}\JM_{n+1}} \\
	\dots \& \dots
	\arrow["{\beta_n}", from=1-1, to=1-2]
	\arrow["{\alpha_n}", from=1-1, to=2-2]
	\arrow["{\beta_n}", from=2-1, to=2-2]
	\arrow["{\alpha_n}", from=2-1, to=3-2]
	\arrow["{\beta_n}", from=3-1, to=3-2]
	\arrow["{\alpha_n}", from=3-1, to=4-2]
\end{tikzcd}\right)
\]
and define $\eta_{n+1}':\mathbb{1}\to P_{n+1}'$ to
be $\eta_n:\mathbb{1}\to
P_{n}$ post-composed with the inclusion of the copy of
$P_n$ in the top right copy of $P_{n}\JM_{n+1}$. Since $P_{n+1}'$
is a twist of a direct sum of copies of $\Cone(\beta_n)$, it obviously kills
all non-identity Bott-Samelsons on the left and right. Now, expanding
$P_n\JM_{n+1}$ to $\Cone(\alpha_n)\to t^{2n+1}P_n$ via rotating triangles, we
see that $\Cone(\eta_{n+1}')$ is of the form
\[\begin{tikzcd}[ampersand replacement=\&,cramped]
	\&\& {(tq)^{-1/2}\mathbb{1}} \\
	{(tq)^{-1/2}qt^{-(2n+1)}P_{n}} \& {t^{-(2n+1)}\Cone(\alpha_n)} \& {P_n} \\
	{(tq)^{-1/2}q^2t^{-2(2n+1)}P_{n}} \& {qt^{-2(2n+1)}P_{n}\JM_{n+1}} \\
	{(tq)^{-1/2}q^3t^{-3(2n+1)}P_{n}} \& {q^2t^{-3(2n+1)}P_{n}\JM_{n+1}} \\
	\dots \& \dots
	\arrow["{\eta_n}"', from=1-3, to=2-3]
	\arrow["{{\beta_n}}", from=2-1, to=2-2]
	\arrow["{{\alpha_n}}", from=2-1, to=3-2]
	\arrow[from=2-2, to=2-3]
	\arrow["{{\beta_n}}", from=3-1, to=3-2]
	\arrow["{{\alpha_n}}", from=3-1, to=4-2]
	\arrow["{{\beta_n}}", from=4-1, to=4-2]
	\arrow["{{\alpha_n}}", from=4-1, to=5-2]
\end{tikzcd}\]
which is a twist of a direct sum of complexes which are all
either $\Cone(\alpha_n)$ or $\Cone(\eta_n)$ and therefore consists only of
non-identity Soergel bimodules. Thus, we have homotopy equivalences
\[\mathbf{P}_{\mathbb{B}_{n+1}}
	\xrightarrow{\id_{\mathbf{P}_{\mathbb{B}_{n+1}}} \otimes
	\eta_{n+1}'} \mathbf{P}_{\mathbb{B}_{n+1}}\otimes P_{n+1}'
	\xleftarrow{\eta_{\mathbb{B}_{n+1}} \otimes \id_{P_{n+1}'}} P_{n+1}'
\]
since the cone of each morphism will be contractible, and so
$P_{n+1}'$ is the unital idempotent from
Theorem~\ref{thm:sbimidempotent}. This proves
Proposition~\ref{prop:symexistence} between $P_n$ and
$\mathbf{P}_{n+1}$.
\subsubsection{Parity of $\End(P_{n+1}')$}
\label{ssec:endparity}
Given the homotopy equivalence $\End(P_{n+1}')\xrightarrow{{\eta_{n+1}'}^*}
\Hom(\mathbb{1},P_{n+1}')$, we can apply the $\Hom(\mathbb{1},-)$ functor to
$P_{n+1}'$ to prove the part of the inductive hypothesis necessary for
Lemma~\ref{lem:homspacesurjection}.
\begin{lemma}
	\label{lem:homspaceparity}
	The complex $\End(P_{n+1}')\simeq \Hom(\mathbb{1},P_{n+1}')$ has homology
	concentrated in cohomological degrees of the same parity as $\End(P_n)$.
\end{lemma}
\begin{proof}
	Take the complex $\Hom(\mathbb{1},\Cone(\beta_n))$, which, by following
	the proof for Corollary~\ref{cor:homcomplex}, will be of the form
\begin{align*}
	\biggl(
		\Hom(\mathbb{1},
		t^{-1/2}q^{1/2}P_n)
		&\xrightarrow{\Hom(\mathbb{1},\id_{P_n}\otimes
	\q[scale=0.4]{\qunit[color=black]})} \Hom (\mathbb{1},q^{1/2}P_nB_{n+1})
	\biggr)\\
	     &{\bigg\downarrow}\\
	     \biggl( \Hom(\mathbb{1},t^{(4n+1)/2}P_nB_{n+1})
	     & \xrightarrow{\Hom(\mathbb{1},\id_{P_n} \otimes
	     \q[scale=0.4]{\qcounit[color=black]})}
     \Hom(P_n,t^{2n+1}P_n) \biggr)
\end{align*}
but since cochain maps $\mathbb{1} \to P_n B_{n+1}$ are in bijection with
cochain maps $\mathbb{1} \to P_n$ via tensoring with a unit, the top row
is an isomorphism and therefore contractible. Additionally,
as the map in the bottom row is an injection, cocycles in
$\Hom(\mathbb{1},t^{(4n+1)/2}P_nB_{n+1})$ cannot contribute to the cohomology
of the complex.
\\ \\
Thus, the complex $\Hom(\mathbb{1},P_{n+1}')$ is of the form
\[\begin{adjustbox}{max width=\textwidth}
	\begin{tikzcd}[ampersand replacement=\&,cramped,column sep=tiny]
	\&\& D \& {\Hom(\mathbb{1},t^{-1/2}P_nB_{n+1})} \& {\Hom(\mathbb{1},P_n)} \\
	{qt^{-(2n+1)}D} \& {\Hom(\mathbb{1},qt^{-(4n+3)}P_nB_{n+1})} \& {\Hom(\mathbb{1},qt^{-(2n+1)}P)} \\
	\dots
	\arrow[from=1-3, to=1-4]
	\arrow[from=1-3, to=2-3]
	\arrow[from=1-4, to=1-5]
	\arrow[from=2-1, to=2-2]
	\arrow[from=2-1, to=3-1]
	\arrow[from=2-2, to=2-3]
\end{tikzcd}
\end{adjustbox}\]
where $D$ is a contractible complex, and due to the directions of the arrows,
the arrows that are injections between the remaining hom complexes are
unchanged by contracting it. Thus, the only terms contributing to the homology
of $\Hom(\mathbb{1},P_{n+1}')$ live in even shifts of the complex
$\Hom(\mathbb{1},P_n)$.
\end{proof}
\noindent This completes the proof of Theorem~\ref{thm:homspaceparity} for
$\mathbf{P}_{n+1}$, and fully completes our induction.
\section{The ring structure of
$H^\bullet \End(\mathbf{P}_{n+1})$}
Our recursion, as well as the proven behavior of the partial trace from
Lemma~\ref{lem:partialtrace} allows us to strengthen
Theorem~\ref{thm:homspaceparity} into the following statement:
\begin{thm}
	\label{thm:ringstructure}
	For $n\geq 1$,
	the ring $H^\bullet\End(\mathbf{P}_{n})$ is isomorphic to
	the bigraded polynomial ring
	\[\mathbb{Q}[u_0,\dots, u_{n-1}]\]
	where $\deg(u_i)=qt^{-(2i+1)}$.
\end{thm}
\begin{proof}
	We proceed by induction. The base case of $n=1$, being identical
	to type $\mathbb{A}$, is already proven.
	\\ \\
	Due to the homotopy equivalence proven in
	Proposition~\ref{prop:symexistence}, we have an isomorphism of
	bigraded rings
	\[H^\bullet\End(\mathbf{P}_{n+1})
	\simeq H^\bullet\End(P_{n+1}').\]
	As $P_{n+1}'$ kills non-identity Soergel bimodules on either side, 
	it absorbs $P_n=I_n(\mathbf{P}_n)$ from the left via the homotopy equivalence
	\[\mathbb{1}\otimes P_{n+1}'\xrightarrow{I_n(\eta_n)\otimes
	\id_{P_{n+1}'}} P_n \otimes P_{n+1}'\]
	giving another isomorphism
	\[H^\bullet \End( P_n\otimes P_{n+1}') \simeq
	H^\bullet\End(P_{n+1}').\]
	Now, the $dg$-ring
	$\End(P_{n}) \simeq \End(\mathbf{P_n})\otimes \mathbb{Q}[x_{n+1}]$
	acts on $\End( P_n\otimes P_{n+1}')$ via post-composition
	(that is, acting by
	$\End(P_n)\otimes \id_{P_{n+1}'}$). By the
	properties of unital idempotents, we have a homotopy equivalence
	\[\End(P_n\otimes P_{n+1}')\xrightarrow{-\circ
	(I_n(\eta_n)\otimes \eta_{n+1}')} \Hom(\mathbb{1},P_n
	\otimes  P_{n+1}')\]
	which is moreover a morphism of $dg$-modules
	of $\End(P_n \otimes P_{n+1}')$, as it is given by pre-composition
	while the module action is given by post-composition. Now, since
	$\pTr_n$ is an
	$\SBim_{\mathbb{B}_n}$-module functor, we have a natural isomorphism
	\[\Hom(\mathbb{1}, P_n\otimes P_{n+1}')\simeq
	\Hom(\mathbb{1}, \mathbf{P}_n\otimes \pTr_n(P_{n+1}')).\]
	The naturality of the isomorphism means that for any
	$g\in \End(P_n)$, the below square commutes:
\[\begin{tikzcd}[ampersand replacement=\&,column sep=huge]
	{\Hom(\mathbb{1},P_n\otimes P_{n+1}')} \& \& 
{\Hom(\mathbb{1},P_n \otimes P_{n+1}')} \\
	{\Hom(\mathbb{1},\mathbf{P}_n\otimes \pTr_n(P_{n+1}'))} \& \&
	{\Hom(\mathbb{1},\mathbf{P}_n\otimes \pTr_n(P_{n+1}'))}
	\arrow["{(g\otimes \id_{P_{n+1}'})\circ -}", from=1-1, to=1-3]
	\arrow[from=1-1, to=2-1]
	\arrow[from=1-3, to=2-3]
	\arrow["{(\pTr_n(g) \otimes \id_{\pTr_n(P_{n+1}')})\circ -}", from=2-1, to=2-3]
\end{tikzcd}\]
making the isomorphism afforded by the adjunction one of
$\End(P_n)$-modules. As proven in~\ref{ssec:endparity},
$\Hom(\mathbb{1},P_{n+1}')$ is homotopy equivalent to
\[\left(
		\Hom(\mathbb{1},t^{-1/2}P_nB_{n+1})
	     \xrightarrow{\Hom(\mathbb{1},\id_{P_n} \otimes
	     \q[scale=0.4]{\qcounit[color=black]})}
	     \Hom(\mathbb{1},P_n)
\right)\otimes \mathbb{Q}[u_{n}]\]
where $u_n$ is the obvious ``row shift'' homomorphism of degree $qt^{-(2n+1)}$.
After applying the partial trace, we observe that this becomes
\[\left(
\begin{tikzcd}[ampersand replacement=\&,cramped]
	\& {\Hom(\mathbb{1},\mathbf{P}_n)} \\
	{\Hom(\mathbb{1},t^{-1/2}q^{1/2}\mathbf{P}_n)} \& {\Hom(1,q\mathbf{P}_n)} \\
	{\Hom(\mathbb{1},t^{-1/2}q^{3/2}\mathbf{P}_n)} \& {\Hom(1,q^2\mathbf{P}_n)} \\
	\dots \& \dots
	\arrow[from=2-1, to=2-2]
	\arrow[from=3-1, to=3-2]
\end{tikzcd}
\right)\otimes \mathbb{Q}[u_{n}]\]
where the arrows are all identities and thus contract, giving a homotopy
equivalence
\[
	\Hom(\mathbb{1},\mathbf{P}_n\otimes \pTr_n(P_{n+1}'))
	\simeq \Hom(\mathbb{1},\mathbf{P}_n) \otimes \mathbb{Q}[u_n]
\]
of $\End(\mathbf{P}_n)$-modules. This proves that
$\End(\mathbf{P}_{n+1})\simeq \End(\mathbf{P}_n)\otimes \mathbb{Q}[u_n]$,
where $u_n$, being a row shift homomorphism, commutes with
all endomorphisms of $\mathbf{P}_n$ and has degree $qt^{-(2n+1)}$.
\end{proof}
\section{The half-$\JM$ in type $\mathbb{B}$}
From the proof of Lemma~\ref{lem:pjmsimplified}, we know that the complex
$P_nF(s_{n+1}\dots s_{2}s_1s_2\dots s_{n+1})$
($P_n$ tensored with half of the generalized Jucys-Murphy braid) will
be homotopy equivalent to
\[\begin{tikzcd}[ampersand replacement=\&,cramped]
	{P_nB_{n+1}\dots B_{212}\dots B_{n+1}} \& \dots \& {t^{(n-1)/2}P_nB_{n+1}\dots B_{212}} \\
	{t^{n/2}P_nB_{n+1}\dots B_{2}B_1} \& \dots \& {t^{(2n+1)/2}P_n}
	\arrow[from=1-1, to=1-2]
	\arrow[from=1-2, to=1-3]
	\arrow[from=1-3, to=2-1]
	\arrow[from=2-1, to=2-2]
	\arrow[from=2-2, to=2-3]
\end{tikzcd}\]
which, after applying the $\Hom(P_n,-)$ functor, will reduce
(via rigid monoidal category properties and Lemma~\ref{lem:homspacezerolightleaves})
to
\[
 \Hom(P_n,t^{n}P_nB_{n+1})
	     \xrightarrow{\Hom(P_n,\id_{P_n} \otimes
	     \q[scale=0.4]{\qcounit[color=black]})}
     \Hom(P_n,t^{(2n+1)/2}P_n)
\]
which is an odd shift of the complex $C$ from the proof of
Lemma~\ref{lem:homspacesurjection}. Unlike in type $\mathbb{A}$,
where tensoring with the
half twist will yield a complex with an odd or even shift relative to $P$,
this complex's shift is always odd, and the homology of $C$ is an odd shift of
$H^\bullet(\End(P_n))/\left<\alpha_{n+1}\right>$, which is simply $\End(P_n)$ in
$\Ch^-(\SBim_{\mathbb{B}_{n}})$ instead of its induced version in
type $\mathbb{B}_{n+1}$.
\section{Further work}
The result of this paper opens up avenues
towards proving categorical diagonalization
statements in type $\mathbb{B}$, in particular proving the diagonalizability of
full-twist Rouquier complexes.
\bibliographystyle{abbrv}
\bibliography{main}

\begin{thebibliography}{1}

\bibitem{abel_categorified_2016}
M.~Abel and M.~Hogancamp.
\newblock Categorified {Young} symmetrizers and stable homology of torus links {II}, Feb. 2016.
\newblock arXiv:1510.05330 [math.QA].

\bibitem{cataltepe_partial_2026}
N.~Cataltepe.
\newblock The partial trace functor on {Soergel} bimodules in type $\mathbb{B}$, Nov. 2026.
\newblock arXiv:2611.25072 [math.QA].

\bibitem{elias_categorical_2017}
B.~Elias and M.~Hogancamp.
\newblock Categorical diagonalization, July 2017.
\newblock arXiv:1707.04349 [math].

\bibitem{elias_computation_2019}
B.~Elias and M.~Hogancamp.
\newblock On the computation of torus link homology.
\newblock {\em Compositio Mathematica}, 155(1):164--205, Jan. 2019.
\newblock arXiv:1603.00407 [math].

\bibitem{elias_categorical_2025}
B.~Elias and M.~Hogancamp.
\newblock Categorical diagonalization of full twists, May 2025.
\newblock arXiv:1801.00191 [math].

\bibitem{elias_introduction_2020}
B.~Elias, S.~Makisumi, U.~Thiel, and G.~Williamson.
\newblock {\em Introduction to {Soergel} {Bimodules}}, volume~5 of {\em {RSME} {Springer} {Series}}.
\newblock Springer International Publishing, Cham, 2020.

\bibitem{hogancamp_idempotents_2017}
M.~Hogancamp.
\newblock Idempotents in triangulated monoidal categories, Mar. 2017.
\newblock arXiv:1703.01001 [math].

\bibitem{hogancamp_categorified_2018}
M.~Hogancamp.
\newblock Categorified {Young} symmetrizers and stable homology of torus links.
\newblock {\em Geometry \& Topology}, 22(5):2943--3002, June 2018.
\newblock arXiv:1505.08148 [math].

\bibitem{hogancamp_constructing_2020}
M.~Hogancamp.
\newblock Constructing categorical idempotents, Feb. 2020.
\newblock arXiv:2002.08905 [math.AT].

\end{thebibliography}
\end{document}